\documentclass[reqno]{amsart}
\usepackage{amssymb,latexsym}
\usepackage{amsfonts,mathrsfs}
\usepackage[margin=1.4in]{geometry}

\newtheorem{theorem}{Theorem}[section]
\newtheorem{lemma}[theorem]{Lemma}
\newtheorem{proposition}[theorem]{Proposition}

\newcommand{\ba}{\begin{array}}
\newcommand{\ea}{\end{array}}

\def \qed{\cqfd}

\newcommand{\Vol}{\mathrm{Vol}}
\newcommand{\rank}{\mathrm{rank}}
\newcommand{\Id}{\mathrm{Id}}

\def\qed{\vbox{\hrule
\hbox{\vrule\hbox to 5pt{\vbox to 8pt{\vfil}\hfil}\vrule}\hrule}}

\newcommand{\beg}{\begin{eqnarray*}}
\newcommand{\begn}{\begin{eqnarray}}
\newcommand{\en}{\end{eqnarray*}}
\newcommand{\enn}{\end{eqnarray}}

\newcommand{\tr}{\mbox{\rm tr\,}}

\begin{document}
\title{A note on curvature estimate of the Hermitian-Yang-Mills flow}
\author{Jiayu Li}
\address{School of Mathematical Sciences,\\
University of Science and Technology of China,\\
Hefei, 230026,\\ and AMSS, CAS, Beijing, 100080, P.R. China\\} \email{jiayuli@ustc.edu.cn}
\author{Chuanjing Zhang}
\address{School of Mathematical Sciences\\
University of Science and Technology of China\\
Hefei, 230026,P.R. China\\ } \email{chjzhang@mail.ustc.edu.cn}
\author{Xi Zhang}
\address{School of Mathematical Sciences,\\
University of Science and Technology of China,\\
Hefei, 230026,P.R. China\\ } \email{mathzx@ustc.edu.cn}
\thanks{The authors were supported in part by NSF in
China, No.11625106, No.11571332 and No. 11131007.}

\begin{abstract} In this paper, we study the curvature estimate of the Hermitian-Yang-Mills flow on holomorphic vector bundles. In one simple case, we show that the curvature of the evolved Hermitian metric is uniformly bounded away from the analytic subvariety determined by the Harder-Narasimhan-Seshadri filtration of the holomorphic vector bundle.

\end{abstract}

\maketitle

\section{Introduction}
\setcounter{equation}{0}

\hspace{0.4cm}

Let $(M ,
\omega )$ be a compact K\"ahler manifold and $\mathcal{E}$ be a coherent sheaf on $M$. A torsion-free coherent sheaf $\mathcal{E}$ is said to be
  $\omega $-stable (respectively, $\omega$-semistable) in Mumford's sense, if for every  coherent proper sub-sheaf
$\mathcal{F}\hookrightarrow \mathcal{E}$, it holds:
\begin{equation}
\mu_{\omega} (\mathcal{F})=\frac{\deg_{\omega} (\mathcal{F})}{\rank (\mathcal{F})}< (\leq ) \mu_{\omega} (\mathcal{E})=\frac{\deg_{\omega} (\mathcal{E})}{\rank (\mathcal{E})},
\end{equation}
where $\mu_{\omega} (\mathcal{F})$ is called the $\omega$-slope of $\mathcal{F}$, the $\omega $-degree of $\mathcal{F}$ is defined as follow
\[
\deg_{\omega}(\mathcal{F})=\int_{M} c_{1}(\mathcal{F})\wedge \frac{\omega ^{n-1}}{(n-1)!},
\]
$c_{1}(\mathcal{F})$ is the first Chern class of $\mathcal{F}$.

To an unstable torsion-free coherent sheaf $\mathcal{E}$, one can associate a
filtration by  subsheaves
\begin{equation}
0=\mathcal{E}_{0}\subset \mathcal{E}_{1}\subset \cdots \subset \mathcal{E}_{k}=\mathcal{E} ,
\end{equation}
such that the quotients $\mathcal{Q}_{i}=\mathcal{E}_{i}/\mathcal{E}_{i-1}$ are torsion-free, $\omega$-semi-stable and $\mu_{\omega} (\mathcal{Q}_{i})>\mu_{\omega} (\mathcal{Q}_{i+1})$, which is called the
Harder-Narasimhan filtration (abbr, HN-filtration) of $\mathcal{E}$. The
associated graded sheaf $Gr^{hn}(\mathcal{E}
)=\oplus_{i=1}^{k}\mathcal{Q}_{i}$ is uniquely determined by the isomorphism
class of $\mathcal{E}$ and the K\"ahler class $[\omega]$. Moreover, for every quotient $\mathcal{Q}_{i}$, there is a further filtration by subsheaves
\begin{equation} 0=\mathcal{E}_{i, 0}\subset
\mathcal{E}_{i, 1}\subset \cdots \subset \mathcal{E}_{i, k_{i}}=\mathcal{Q}_{i} ,\end{equation}  such that the quotients $\mathcal{Q}_{i,j}=\mathcal{E}_{i, j}/\mathcal{E}_{i, j-1}$ is
torsion-free and $\omega$-stable, $\mu_{\omega} (\mathcal{Q}_{i, j})=\mu
(\mathcal{Q}_{i})$ for each $j$. This double filtration $\{\mathcal{E}_{i, j}\}$ is called the Harder-Narasimhan-Seshadri filtration (abbr, HNS-filtration) of the sheaf $\mathcal{E}$. The
associated graded sheaf: $ Gr^{hns}(E, \overline{\partial}_{A}, \phi )=\oplus_{i=1}^{k}\oplus_{j=1}^{k_{i}}Q_{i , j}
$ is uniquely determined by the isomorphism class of $\mathcal{E}$ and the K\"ahler class $[\omega]$. The number $\sum_{i=1}^{k}k_{i}-1$ is called the length of the HNS-filtration.

In the following, we denote $\Sigma_{\mathcal{E}}$ the set of singularities where $\mathcal{E}$ is not locally free.  If $\mathcal{E}$ is locally free on the whole $M$, i.e. $\Sigma_{\mathcal{E}}=\varnothing $, there is a holomorphic vector bundle $(E, \overline{\partial }_{E})$ on $M$ such that the sheaf $\mathcal{E}$ is generated by the local holomorphic sections of $(E, \overline{\partial }_{E})$. A locally free coherent sheaf $\mathcal{E}$ can be seen as a holomorphic vector bundle, i.e. $\mathcal{E}=(E, \overline{\partial }_{E})$. A Hermitian metric $H$ on the  holomorphic vector bundle $(E, \overline{\partial }_{E})$ is said to be
$\omega$-Hermitian-Einstein  if it
  satisfies the following Einstein condition on $M $, i.e.
\begin{equation}\label{HE}
\sqrt{-1}\Lambda_{\omega} F_{H}
=\lambda_{E, \omega} \Id_{\mathcal{E}},
\end{equation}
where   $\lambda_{E, \omega}  =\frac{2\pi}{\Vol(M, \omega)} \mu_{\omega} (\mathcal{E})$,
$F_{H}$ is the curvature tensor of Chern connection $D_{H}$ with respect to the Hermitian metric $H$, and $\Lambda_{\omega }$ denotes the contraction with  the K\"ahler metric $\omega $.
The Donaldson-Uhlenbeck-Yau theorem (\cite{NS,D0,D1,UY}) states that the $\omega$-stability of $\mathcal{E}$ implies the existence of  $\omega$-Hermitian-Einstein metric on $\mathcal{E}$. This theorem has several interesting and important generalizations and extensions (\cite{LY, Hi, Si, Br1, BS, GP, BG, HL, AG1, Bi, BT, JZ, LN1, LN2, M, Mo1, Mo}, etc.).

\hspace{0.3cm}

Let $H_{0}$ be  a Hermitian metric  on the complex vector bundle $E$,
$ \mathcal{A}_{H_{0}}$ be the space of connections of $E$
compatible with the metric $H_{0}$, and $ \mathcal{A}^{1,1}_{H_{0}}$
be the space of unitary integrable connections of $E$ (i.e. those whose curvature is of type $(1, 1)$). For any $A_{0}\in \mathcal{A}^{1,1}_{H_{0}}$, $\overline{\partial}_{A_{0}}=D_{A_{0}}^{(0, 1)}$
defines a
holomorphic structure on $E$. We consider the following Yang-Mills flow on the Hermitian vector bundle $(E, H_{0})$ with initial data $A_{0}$,
\begin{equation}\label{YMF}
\left \{\begin{split}
& \frac{\partial A}{\partial t}=-D_{A}^{\ast} F_{A},\\
& A(0)=A_{0}.\\
\end{split}
\right.
\end{equation}
The Yang-Mills flow, as the gradient flow of the Yang-Mills functional, was first suggested by Atiyah-Bott in \cite{AB}.
Donaldson \cite{D1}  proves the global existence of the Yang-Mills flow if the initial data $A_{0}$ is integrable. In fact, Donaldson introduces the following Hermitian-Yang-Mills flow for Hermitian metrics $H(t)$ on the holomorphic bundle $(E, \overline{\partial}_{A_{0}})$ with initial metric $H_{0}$,
\begin{equation}\label{Dflow}
\left \{\begin{split} &H^{-1}\frac{\partial H}{\partial t}=-2 (\sqrt{-1}\Lambda_{\omega }F_{H} -\lambda \Id_{E}),\\
&H(0)=H_{0}, \\
\end{split}
\right.
\end{equation}
and proves that  solution to the above nonlinear heat equation exists for all time. Then Donaldson shows that,  by choosing complex gauge transformations $\sigma (t)$ which satisfy $\sigma(t)^{\ast }\sigma (t)=H_{0}^{-1}H(t)$, $A(t)=\sigma (t)(A_{0})$ is the unique long time solution of the Yang-Mills flow (\ref{YMF}). Furthermore, Donaldson proves the
convergence of the flow at infinity in the case that the
initial holomorphic structure $(E, \overline{\partial}_{A_{0}})$ is stable.

In general case, since the mean curvature tensor $\sqrt{-1}\Lambda_{\omega }F_{A(t)}$ is uniformly bounded along the Yang-Mills flow, the Uhlenbeck's compactness result (\cite{Uh1,Uh2}) implies that for any sequence $A(t_{i})$ along the flow there is a subsequence, modulo gauge transformations, weakly converges  to a Yang-Mills connection $A_{\infty }$ outside a closed subset $\Sigma_{an}$ of Hausdorff complex codimension at least $2$. We call $\Sigma_{an}$ the bubbling set. In \cite{HT}, Hong and Tian have shown that in fact the convergence  can be taken to be $C^{\infty}_{loc}$ on $M\setminus \Sigma_{an}$ and the bubbling set $\Sigma_{an}$ is a holomorphic subvariety. By Bando and Siu's result (\cite{BS}), the holomorphic vector bundle $(E_{\infty}, \overline{\partial}_{A_{\infty}})$ on $M\setminus \Sigma_{an}$ can be extended to the whole $M$ as a reflexive coherent sheaf $\mathcal{E}_{\infty}$. On the other hand, since  $A_{\infty}$ is Yang-Mills, $\mathcal{E}_{\infty}$ has a  holomorphic orthogonal splitting of stable reflexive sheaves with admissible Hermitian-Einstein metrics. In \cite{BS}, Bando and Siu conjecture that this limiting reflexive coherent sheaf $\mathcal{E}_{\infty}$ should be isomorphic to the double dual of the graded object of the Harder-Narasimhan-Seshadri filtration of the initial holomorphic structure $(E, \overline{\partial}_{A_{0}})$. This isomorphism was pointed out by Atiyah-Bott  \cite{AB} first for the  Riemann surface case, and proved  by Daskalopoulos \cite{Da}. Later, this conjecture was proved by Daskalopoulos and Wentworth \cite{DW1} in the case of  K\"ahler surfaces, and Jacob \cite{Ja3} and Sibley \cite{Sib}  for the higher dimension case. This conjecture is also set up in the Higgs bundle case, see references \cite{Wi,LZ1,LZ2}.

In the following, we denote $\Sigma_{alg}$ the singular set of  the associated graded sheaf $Gr^{hns}(E, \overline{\partial}_{A_{0}})$, i.e. $Gr^{hns}(E, \overline{\partial}_{A_{0}})$ is
locally free away from $\Sigma_{alg}$. $\Sigma_{alg}$ is a complex analytic subvariety of complex codimension at least $2$, which we call the algebraic singular set. By the results in \cite{Sib,Ja3}, it is not hard to see that $\Sigma_{alg} \subset \Sigma_{an}$. It is an interesting problem to prove that $\Sigma_{an}\subset \Sigma_{alg}$. This problem was solved by Daskalopoulos and Wentworth \cite{DW2} for $\dim_{\mathbb{C}}M=2$, Sibley and Wentworth (\cite{SW}) for the higher dimensions case. Sibley and Wentworth's method are mostly algebraic, it should be an interesting problem to give a uniform curvature estimate of the Yang-Mills flow  away from the algebraic singular set $\Sigma_{alg}$ by using analytic methods. In this paper, we solve the problem in the case that $(E, \overline{\partial}_{A_{0}})$ is nonsemistable and the Harder-Narasimhan-Seshadri filtration is of length one. In fact we obtain the following theorem.

\medskip

\begin{theorem}\label{UCE}
Let $(E, H_{0})$ be a Hermitian vector bundle over a compact K\"ahler manifold $(M, \omega)$, and $A(t)$ be the solution of the Yang-Mills flow (\ref{YMF}) with initial connection $A_{0}\in \mathcal{A}^{1,1}_{H_{0}}$. If the holomorphic bundle   $\mathcal{E}=(E, \overline{\partial}_{A_{0}})$ is nonsemistable and the Harder-Narasimhan-Seshadri filtration is of length one. Then, for any compact subset $U\subset M\setminus \Sigma_{alg}$ there exists a uniform constant $C_{U}$ such that
\begin{equation}\label{CE}
\sup_{(x, t)\in U\times [0, +\infty )}|F_{A}|_{H_{0}}^{2}(x, t)\leq C_{U}.
\end{equation}
\end{theorem}

\medskip

We believe the theorem holds in general, the proof is more complicated.
We now give an overview of our proof.
Let $H(t)$ be the long-time solution of the Hermitian Yang-Mills flow (\ref{Dflow}), it is well known that
\begin{equation}
|F_{H(t)}|_{H(t)}^{2}=|F_{A}|_{H_{0}}^{2}.
\end{equation}
So, we need only to estimate the curvature tensor $F_{H(t)}$ of Chern connection $D_{H(t)}$ with respect to the evolved Hermitian metric $H(t)$. By the assumption of theorem \ref{UCE}, there exists an exact sequence
\begin{equation}
0\rightarrow \mathcal{S}\rightarrow \mathcal{E} \rightarrow \mathcal{Q} \rightarrow 0
\end{equation}
such that $\mathcal{S}$ and $\mathcal{Q}$ are torsion-free  $\omega$-stable  sheaves, and
\begin{equation}
\mu_{\omega }(\mathcal{S})>\mu_{\omega }(\mathcal{E})>\mu_{\omega }(\mathcal{Q}).
\end{equation}
On $M\setminus \Sigma_{alg}$, $\mathcal{S}$ and $\mathcal{Q}$ can be seen as holomorphic vector bundles, we denote $H_{S}(t)$ and $H_{Q}(t)$ the Hermitian metrics on the subbundle $\mathcal{S}$ and the quotient bundle $\mathcal{Q}$ induced by the evolved metric $H(t)$ on $E$, and $\gamma (t)$ the second fundamental form. We derive the evolution equations for $H_{S}(t)$, $H_{Q}(t)$ and $\gamma (t)$ (see (\ref{S0}), (\ref{Q0}) and (\ref{r0})).  In \cite{Si}, Simpson generalizes the Hermitian-Yang-Mills flow to Higgs bundle case. Under the assumption of stability, by using a result of Uhlenbeck and Yau (\cite{UY}), Simpson obtains a uniform $C^{0}$-estimate on $H(t)$, this implies   uniform higher order estimates including the uniform curvature estimate. When $\mathcal{E}$ is unstable, the $C^{0}$-norm of the evolved metrics $H(t)$ might be unbounded. In our case, we can obtain a uniform local $C^{0}$-bound on the rescaled metrics $\tilde{H}_{S}(t)=e^{2(\lambda_{S}-\lambda_{E})t}H_{S}(t)$ and $\tilde{H}_{Q}(t)=e^{2(\lambda_{Q}-\lambda_{E})t}H_{Q}(t)$ away from the algebraic singular set $\Sigma_{alg}$. This uniform local $C^{0}$-estimates is a key point in the proof of our theorem, where we will use the stabilities of $\mathcal{S}$ and $\mathcal{Q}$ and the property that $\mu_{\omega }(\mathcal{S})>\mu_{\omega }(\mathcal{Q})$, see section $3$ for details. By using the above local $C^{0}$-estimates, we prove that the norms of the second fundamental forms $|\gamma (t)|_{H(t)}$, $|T_{S}(t)|_{H_{S}(t)}$ and $|T_{Q}(t)|_{H_{Q}(t)}$ are uniformly bounded away from $\Sigma_{alg}$. By choosing suitable test functions and using the maximum principle, we obtain a uniform local estimate of $|F_{H(t)}|_{H(t)}$ on $M\setminus \Sigma_{alg}$.
In \cite{CJ}, under some technical assumptions on the growth of norms of the second fundamental forms associated to the HNS filtration, Collins and Jacob  obtain the above uniform curvature estimate.

\medskip

This paper is organized as follows. In Section 2, we recall some basic estimates for the Hermitian-Yang-Mills flow, and derive the evolution equations for the induced metrics $H_{S}(t)$ and $H_{Q}(t)$ and the second fundamental forms $\gamma (t)$.   In section 3, we recall the resolution of the HNS filtration of holomorphic vector bundle, prove the related Donaldson's functionals are uniformly bounded from below, and obtain a uniform local $C^{0}$-bound on the rescaled metrics.  In section 4, we obtain a uniform local estimate for the norms of the second fundamental forms and a uniform local $C^{1}$-estimate for the rescaled metrics. Then we complete the proof of Theorem \ref{UCE} in section 5.

\section{Evolution of the second fundamental form}
\setcounter{equation}{0}

Let $(M, \omega )$ be a K\"ahler manifold which may be noncompact, and $(E, \overline{\partial}_{E})$ be a holomorphic vector bundle on $M$. We suppose that there is an exact sequence of holomorphic vector bundles:
\begin{equation}
0\rightarrow S\rightarrow E\rightarrow Q \rightarrow 0.
\end{equation}
We denote $\overline{\partial }_{S}$ and $\overline{\partial }_{Q}$ the holomorphic structures on the holomorphic subbundle $S$ and the quotient bundle $Q$ induced by $\overline{\partial }_{E}$, $\overline{\partial }_{S\oplus Q}$ the induced holomorphic structure on the direct sum bundle $S\oplus Q$, i.e.
\begin{equation}
\overline{\partial }_{S\oplus Q}=\begin{pmatrix}
\overline{\partial}_{S} &  0 \\
0   & \overline{\partial}_{Q}\\
\end{pmatrix}.
\end{equation}

For every Hermitian metric $H$ on $E$, we have the following bundle isomorphism
\begin{equation}\label{is1}
f_{H}: S\oplus Q \rightarrow E , \qquad (X, [Y])\mapsto i(X) +(\Id_{E}-\pi _{H})(Y),
\end{equation}
 where $X\in S$, $Y\in E $, $i:S\hookrightarrow E$ is the inclusion and $\pi_{H}: E\rightarrow E$ is the orthogonal projection into $S$ with respect to the metric $H$.
In the following, we denote $H_{S}$ and $H_{Q}$ the Hermitian metrics on $S$ and $Q$ induced by the metric $H$ on $E$. By the definition, the pulling back metric is
\begin{equation}
f_{H}^{\ast}(H)=
\begin{pmatrix}
H_{S} &  0 \\
0   & H_{Q}\\
\end{pmatrix}.
\end{equation}
Since $\overline{\partial }_{E}^{2}=0$, the pull back holomorphic structure
\begin{equation}
f_{H}^{\ast }(\overline{\partial }_{E})=f_{H}^{-1} \circ \overline{\partial }_{E}\circ f_{H}
\end{equation}
is also a holomorphic structure
on $S\oplus Q$.
Recall that $S$ is a holomorphic subbundle of $(E, \overline{\partial }_{E})$, for any $e\in S$ and $[Y]\in Q$, we have
\begin{equation}
f_{H}^{\ast }(\overline{\partial }_{E})(e, [Y])= (\overline{\partial }_{S}e +\gamma ([Y]), \overline{\partial }_{Q}[Y]),
\end{equation}
where $\gamma ([Y])=-(\overline{\partial }_{E} \pi_{H})(Y)$.
So we have the following expression
\begin{equation}\label{Se1}
f_{H}^{\ast }(\overline{\partial }_{E})-\overline{\partial }_{S\oplus Q}=
\begin{pmatrix}
0 &  \gamma  \\
0   & 0\\
\end{pmatrix},
\end{equation}
where $\gamma (t)\in \Omega ^{0, 1}(Hom (Q, S))$  will be called  the  second fundamental form. Furthermore,
\begin{equation}
\overline{\partial }_{S\otimes Q^{\ast}}\gamma =\overline{\partial }_{S}\circ \gamma +\gamma \circ \overline{\partial }_{Q}=0.
\end{equation}
Let $D_{(\overline{\partial}_{E}, H) }$ be the Chern connection determined by the holomorphic structure $\overline{\partial }_{E}$ and the metric $H$, it is easy to see that the pulling back connection $f_{H}^{\ast}(D_{(\overline{\partial}_{E}, H) })$ is just the Chern connection on the holomorphic bundle $(S\oplus Q, f_{H}^{\ast }(\overline{\partial }_{E}))$ with respect to the metric $f_{H}^{\ast}(H)$, so
\begin{equation}
f_{H}^{\ast }(\partial_{H})=\left (
\begin{matrix}
\partial_{H_{S}} &  0 \\
-\gamma (t)^{\ast H}   & \partial_{H_{Q}}\\
\end{matrix}
\right ),
\end{equation}
where $\partial_{H}=D^{1, 0}_{(\overline{\partial}_{E}, H)}$ , $\partial_{H_{S}}=D^{1, 0}_{(\overline{\partial}_{S}, H_{S})}$, $\partial_{H_{Q}}=D^{1, 0}_{(\overline{\partial}_{Q}, H_{Q})}$ and $\gamma (t)^{\ast H}\in \Omega ^{1, 0}(S^{\ast}\otimes Q)$ is the adjoint of $\gamma $ with respect to metrics $H_{S}$ and $H_{Q}$.
We also have the pulling back curvature form, i.e. the following Gauss-Codazzi equation
 \begin{equation}\label{GC}
 f_{H}^{\ast }(F_{H})=
\begin{pmatrix}
F_{H_{S}}-\gamma \wedge \gamma^{\ast H}  & \partial_{H_{S}} \gamma +\gamma \partial_{H_{Q}} \\
-\overline{\partial }_{Q}\gamma^{\ast H}-\gamma^{\ast H}\overline{\partial }_{S}   & F_{H_{Q}}-\gamma^{\ast H} \wedge \gamma \\
\end{pmatrix},
\end{equation}
where $F_{H_{S}}$ and $F_{H_{Q}}$ are the curvature forms of the Chern connections $D_{(\overline{\partial}_{S}, H_{S})}$ and $D_{(\overline{\partial}_{Q}, H_{Q})}$.

Let $H(t)$ be the solution of the Hermitian-Yang-Mills flow (\ref{Dflow}) on the holomorphic bundle $(E, \overline{\partial }_{E})$ with initial metric $H_{0}$.
Now, we will split the Hermitian-Yang-Mills flow (\ref{Dflow}) to the subbundle $S$ and  the quotient bundle $Q$.
Since $f_{H(t)}|_{S}=i $,  we have
\begin{equation}
\frac{\partial f_{H(t)}}{\partial t}(e)=\frac{\partial }{\partial t} (f_{H(t)}(e))=0,
\end{equation}
\begin{equation}
\frac{\partial f_{H(t)}}{\partial t}([Y])=\frac{\partial }{\partial t} ((\Id -\pi_{H(t)})Y)\in S,
\end{equation}
for every  $e\in S$ and $Y\in E$.
So, we have the following expression
\begin{equation}
f_{H(t)}^{-1 }\frac{\partial f_{H(t)}}{\partial t} (e_{1}, \cdots, e_{s}, \xi_{1}, \cdots , \xi_{q})=(e_{1}, \cdots, e_{s}, \xi_{1}, \cdots , \xi_{q})\left(\begin{matrix}
0 & \chi (t) \\
0   & 0\\
\end{matrix}
\right),
\end{equation}
where $\{e_{j}\}_{j=1}^{s}$ ($\{\xi_{\alpha }\}_{\alpha =1}^{q}$) is a local basis of the bundle $S$ (resp. $Q$), and $\chi (t) \in \Gamma (S\otimes Q^{\ast})$.
For simplicity, we denote $\bar{H}(t)=f_{H(t)}^{\ast}(H(t))$. For every $X, Y \in S\oplus Q$, we have
\begin{equation}
\begin{split}
\langle\bar{H}^{-1}\frac{\partial \bar{H}}{\partial t}(X), Y\rangle_{\bar{H}}=&\frac{\partial }{\partial t}\bar{H} (X, Y)= \frac{\partial }{\partial t} \langle X, Y\rangle_{\bar{H}}\\
=& \frac{\partial }{\partial t} \langle f_{H}(X), f_{H}(Y)\rangle_{H}\\
=& \frac{\partial H}{\partial t} (f_{H}(X), f_{H}(Y))+\langle\frac{\partial f_{H}}{\partial t}(X), f_{H}(Y)\rangle_{H}\\
&+\langle f_{H}(X), \frac{\partial f_{H}}{\partial t}(Y)\rangle_{H}\\
=&\langle H^{-1}\frac{\partial H}{\partial t} (f_{H}(X)), f_{H}(Y)\rangle_{H}+\langle f_{H}^{-1}\frac{\partial f_{H}}{\partial t}(X), Y\rangle_{\bar{H}}\\
&+\langle X, f_{H}^{-1}\frac{\partial f_{H}}{\partial t}(Y)\rangle_{\bar{H}}\\
=&\langle f_{H}^{\ast}(H^{-1}\frac{\partial H}{\partial t}) (X), Y\rangle_{\bar{H}}+\langle f_{H}^{-1}\frac{\partial f_{H}}{\partial t}(X), Y\rangle_{\bar{H}}\\
&+\langle X, f_{H}^{-1}\frac{\partial f_{H}}{\partial t}(Y)\rangle_{\bar{H}},
\end{split}
\end{equation} and then
\begin{equation}
f_{H}^{\ast }(H^{-1}\frac{\partial H}{\partial t})=\left(\begin{matrix}
H_{S}^{-1}\frac{\partial H_{S}}{\partial t} & 0 \\
0   & H_{Q}^{-1}\frac{\partial H_{Q}}{\partial t}\\
\end{matrix}
\right)-f_{H}^{-1 }\frac{\partial f_{H}}{\partial t} -(f_{H}^{-1 }\frac{\partial f_{H}}{\partial t})^{\ast \bar{H}}.
\end{equation}

By (\ref{Dflow})  and the Gauss-Codazzi equation (\ref{GC}), we have
\begin{equation}\label{S0}
H_{S}^{-1}\frac{\partial H_{S}}{\partial t}=-2(\sqrt{-1}\Lambda_{\omega} (F_{H_{S}}-\gamma \wedge \gamma^{\ast}) -\lambda \Id),
\end{equation}
\begin{equation}\label{Q0}
H_{Q}^{-1}\frac{\partial H_{Q}}{\partial t}=-2(\sqrt{-1}\Lambda_{\omega} (F_{H_{Q}}-\gamma^{\ast} \wedge \gamma ) -\lambda \Id),
\end{equation}
\begin{equation}\label{f0}
f_{H}^{-1 }\frac{\partial f_{H}}{\partial t} =\left(\begin{matrix}
0 & 2\sqrt{-1}\Lambda_{\omega }(\partial_{H_{S}}\gamma + \gamma \partial_{H_{Q}}) \\
0   & 0\\
\end{matrix}
\right).
\end{equation}

Now, we consider the evolution of the second fundamental form $\gamma (t)$, i.e. we have:

\medskip

\begin{lemma}  Let $H(t)$ be the solution of the heat flow (\ref{Dflow})  with initial metric $H_{0}$ and $\gamma (t)$ be the second fundamental form defined by the formula (\ref{Se1}), then we have:
\begin{equation}\label{r0}
\frac{\partial }{\partial t} \gamma =2\overline{\partial }_{S\otimes Q^{\ast}}(\sqrt{-1}\Lambda_{\omega }(\partial_{H_{S}}\gamma + \gamma \partial_{H_{Q}})).
\end{equation}
\end{lemma}

\medskip

{\bf Proof.}
For simplicity, we denote
\begin{equation}
\overline{\partial }f_{H}=\overline{\partial }_{E} \circ f_{H}- f_{H} \circ \overline{\partial }_{S\oplus Q}
\end{equation}
and then
\begin{equation}
f_{H}^{-1}\overline{\partial }f_{H}= f_{H}^{\ast}(\overline{\partial }_{E})- \overline{\partial }_{S\oplus Q}=
\left(\begin{matrix}
0 & \gamma \\
0   & 0\\
\end{matrix}
\right).
\end{equation}
Taking the derivative of the above equation with respect to $t$, then
\begin{equation}
\begin{split}
\left(\begin{matrix}
0 & \frac{\partial \gamma }{\partial t} \\
0   & 0\\
\end{matrix}
\right)=& \frac{\partial }{\partial t}(f_{H}^{-1}\overline{\partial }f_{H})\\
=& -f_{H}^{-1}\frac{\partial f_{H}}{\partial t}f_{H}^{-1}\overline{\partial }f_{H}+ f_{H}^{-1}\frac{\partial }{\partial t}(\overline{\partial }f_{H})\\
=& -f_{H}^{-1}\frac{\partial f_{H}}{\partial t}f_{H}^{-1}\overline{\partial }f_{H}+ f_{H}^{-1}\circ \overline{\partial }_{E}\circ f_{H}\circ f_{H}^{-1}\frac{\partial f_{H}}{\partial t}\\
&- f_{H}^{-1}\frac{\partial f_{H}}{\partial t}\circ \overline{\partial }_{S\oplus Q}\\
=& -\left(\begin{matrix}
0 & 2\sqrt{-1}\Lambda_{\omega }(\partial_{H_{S}}\gamma + \gamma \partial_{H_{Q}}) \\
0   & 0\\
\end{matrix}
\right)\left(\begin{matrix}
0 & \gamma \\
0   & 0\\
\end{matrix}
\right) \\
&+\left(\begin{matrix}
\overline{\partial}_{S} &  \gamma \\
0   & \overline{\partial}_{Q}\\
\end{matrix}
\right )\left(\begin{matrix}
0 & 2\sqrt{-1}\Lambda_{\omega }(\partial_{H_{S}}\gamma + \gamma \partial_{H_{Q}}) \\
0   & 0\\
\end{matrix}
\right)\\
&-\left(\begin{matrix}
0 & 2\sqrt{-1}\Lambda_{\omega }(\partial_{H_{S}}\gamma + \gamma \partial_{H_{Q}}) \\
0   & 0\\
\end{matrix}
\right)
\left(\begin{matrix}
\overline{\partial}_{S} &  0 \\
0   & \overline{\partial}_{Q}\\
\end{matrix}
\right )\\
=& \left(\begin{matrix}
0 & 2\overline{\partial }_{S\otimes Q^{\ast}}(\sqrt{-1}\Lambda_{\omega }(\partial_{H_{S}}\gamma + \gamma \partial_{H_{Q}})) \\
0   & 0\\
\end{matrix}
\right),
\end{split}
\end{equation}
i.e. we have the formula (\ref{r0}).

\hfill $\Box$ \\

\medskip

\begin{lemma}  Let $f_{H(t)}$ be the bundle isomorphism defined in (\ref{is1}), then we have
\begin{equation}\label{is2}
f_{H_{0}}^{-1} f_{H(t)}=\left(\begin{matrix}
\Id_{S} & G(t) \\
0   & \Id_{Q}\\
\end{matrix}
\right),
\end{equation}
where $G (t)\in \Gamma (S\otimes Q^{\ast})$ and $G(0)=0$. Furthermore,
\begin{equation}\label{G0}
\frac{\partial }{\partial t} G =2\sqrt{-1}\Lambda_{\omega }(\partial_{H_{S}}\gamma + \gamma \partial_{H_{Q}}) ,
\end{equation}
and
\begin{equation}\label{G1}
\overline{\partial }_{S\otimes Q^{\ast}}G=\gamma -\gamma_{0}.
\end{equation}

 \end{lemma}

\medskip

{\bf Proof. }
By the definition, we have
\begin{equation}
f_{H_0}^{-1} f_{H}(e)=f_{H_0}^{-1} (e)=e,
\end{equation}
and
\begin{equation}
f_{H_0}^{-1} f_{H}([Y])=f_{H_0}^{-1} ((\Id -\pi_{H})Y)=[Y]+(\pi_{H_{0}}-\pi_{H})Y,
\end{equation}
for every  $e\in S$ and $Y\in E$.
So we have the expression (\ref{is2}).

 Taking the derivative of the equation (\ref{is2}) with respect to $t$ and using the formula (\ref{f0}), we get
\begin{equation}
\begin{split}
\frac{\partial }{\partial t}(f_{0}^{-1} f_{H})=&f_{0}^{-1}f_{H} \circ f_{H}^{-1}\frac{\partial f_{H} }{\partial t}\\
=& \left(\begin{matrix}
\Id_{S} & G \\
0   & \Id_{Q}\\
\end{matrix}
\right)\left(\begin{matrix}
0 & 2\sqrt{-1}\Lambda_{\omega }(\partial_{H_{S}}\gamma + \gamma \partial_{H_{Q}}) \\
0   & 0\\
\end{matrix}
\right)\\
=&\left(\begin{matrix}
0 & 2\sqrt{-1}\Lambda_{\omega }(\partial_{H_{S}}\gamma + \gamma \partial_{H_{Q}}) \\
0   & 0\\
\end{matrix}
\right),
\end{split}
\end{equation}
i.e.  we obtain the formula (\ref{G0}).

On the other hand, we have
\begin{equation}
\begin{split}
\left(\begin{matrix}
0 & \overline{\partial }_{S\otimes Q^{\ast }}G \\
0   & 0\\
\end{matrix}
\right)=&
\left(\begin{matrix}
0 & \overline{\partial }_{S}\circ G-G\circ \overline{\partial }_{Q} \\
0   & 0\\
\end{matrix}
\right)\\
=& \left(\begin{matrix}
\overline{\partial}_{S} &  0 \\
0   & \overline{\partial}_{Q}\\
\end{matrix}
\right )\left(\begin{matrix}
\Id_{S} & G \\
0   & \Id_{Q}\\
\end{matrix}
\right)- \left(\begin{matrix}
\Id_{S} & G \\
0   & \Id_{Q}\\
\end{matrix}
\right)\left(\begin{matrix}
\overline{\partial}_{S} &  0 \\
0   & \overline{\partial}_{Q}\\
\end{matrix}
\right )\\
=& \overline{\partial }_{S\oplus Q}(f_{0}^{-1} f_{H})\\
=&\overline{\partial }_{S\oplus Q} \circ  f_{0}^{-1} f_{H} - f_{0}^{-1} f_{H} \circ \overline{\partial }_{S\oplus Q}\\
=&\overline{\partial }_{S\oplus Q} \circ  f_{0}^{-1} f_{H}-f_{0}^{-1}\circ \overline{\partial }_{E}\circ  f_{H}\\ &+f_{0}^{-1}\circ \overline{\partial }_{E}\circ  f_{H} - f_{0}^{-1} f_{H} \circ \overline{\partial }_{S\oplus Q}\\
=&\overline{\partial }_{S\oplus Q} \circ  f_{0}^{-1} f_{H}-f_{0}^{\ast}( \overline{\partial }_{E})\circ  f_{0}^{-1}f_{H}\\
&+f_{0}^{-1}f_{H}\circ f_{H}^{\ast}(\overline{\partial }_{E}) - f_{0}^{-1} f_{H} \circ \overline{\partial }_{S\oplus Q}\\
=&\{\overline{\partial }_{S\oplus Q} -f_{0}^{\ast}( \overline{\partial }_{E})\}\circ  f_{0}^{-1}f_{H}+f_{0}^{-1}f_{H}\circ \{f_{H}^{\ast}(\overline{\partial }_{E}) -  \overline{\partial }_{S\oplus Q}\}\\
=&- \left(\begin{matrix}
0 &  \gamma_{0} \\
0   & 0\\
\end{matrix}
\right )\left(\begin{matrix}
\Id_{S} & G \\
0   & \Id_{Q}\\
\end{matrix}
\right)+ \left(\begin{matrix}
\Id_{S} & G \\
0   & \Id_{Q}\\
\end{matrix}
\right)\left(\begin{matrix}
0 &  \gamma \\
0   & 0\\
\end{matrix}
\right )\\
=& \left(\begin{matrix}
0 & \gamma -\gamma_{0}\\
0   & 0\\
\end{matrix}
\right),
\end{split}
\end{equation}
so we obtain the formula (\ref{G1}).

\hfill $\Box$ \\

In the following, we consider the parabolic inequalities  for  $|\gamma (t)|_{H(t)}^{2}= -\sqrt{-1}\Lambda_{\omega } \tr( \gamma \wedge \gamma ^{\ast H(t)})$ and $|G(t)|_{H(t)}^{2}= \tr( G \circ G ^{\ast H(t)})$ , which will be needed in the next section.
By direct calculations, we get
\begin{equation}
\begin{split}
\Delta |\gamma (t)|_{H(t)}^{2}=& 2g^{k\bar{l}}\partial_{k}\overline{\partial }_{l} |\gamma (t)|_{H(t)}^{2}\\
=& 2g^{k\bar{l}}\{ \langle\nabla _{\partial_{k}}^{H(t)}\nabla _{\overline{\partial}_{l}}^{H(t)}\gamma , \gamma \rangle_{H(t)}+\langle\gamma , \nabla _{\overline{\partial}_{k}}^{H(t)}\nabla _{\partial_{l}}^{H(t)}\gamma  \rangle_{H(t)}\\
& +\langle\nabla _{\partial_{k}}^{H(t)}\gamma , \nabla _{\partial_{l}}^{H(t)}\gamma \rangle_{H(t)} +\langle\nabla _{\overline{\partial}_{l}}^{H(t)}\gamma , \nabla _{\overline{\partial}_{l}}^{H(t)}\gamma \rangle_{H(t)}\}\\
=& 2Re \{g^{k\bar{l}} \langle\nabla _{\partial_{k}}^{H(t)}\nabla _{\overline{\partial}_{l}}^{H(t)}\gamma +\nabla _{\overline{\partial}_{l}}^{H(t)}\nabla _{\partial_{k}}^{H(t)}\gamma , \gamma \rangle_{H(t)}\}\\
&+2|\nabla ^{H(t)}\gamma |_{H(t)}^{2},\\
\end{split}
\end{equation}
\begin{equation}
\begin{split}
\nabla _{\overline{\partial}_{l}}^{H(t)}\gamma =&\nabla _{\overline{\partial}_{l}}^{H(t)}(\gamma_{\bar{j}} d\overline{z}^{j})\\
=& (D _{\overline{\partial}_{l}}^{H(t)}\gamma_{\bar{j}}) d\overline{z}^{j} +\gamma_{\bar{j}}\nabla _{\overline{\partial}_{l}}d\overline{z}^{j}, \\
\end{split}
\end{equation}

\begin{equation}
\begin{split}
\nabla _{\partial _{k}}^{H(t)}\nabla _{\overline{\partial}_{l}}^{H(t)}\gamma  =& (D _{\partial _{k}}^{H(t)}D _{\overline{\partial}_{l}}^{H(t)}\gamma_{\bar{j}}) d\overline{z}^{j} +\gamma_{\bar{j}}(\nabla _{\partial _{k}}\nabla _{\overline{\partial}_{l}}d\overline{z}^{j})\\
&+(D _{\partial _{k}}^{H(t)}\gamma_{\bar{j}})\nabla _{\overline{\partial}_{l}}d\overline{z}^{j}, \\
\end{split}
\end{equation}
and
\begin{equation}
\begin{split}
 &D _{\partial _{k}}^{H(t)}D _{\overline{\partial}_{l}}^{H(t)}\gamma_{\bar{j}}-D _{\overline{\partial}_{l}}^{H(t)}D _{\partial _{k}}^{H(t)}\gamma_{\bar{j}}\\
 =& F_{H_{S}}(\partial _{k} , \overline{\partial}_{l} )\gamma_{\bar{j}} -\gamma_{\bar{j}}F_{H_{Q}}(\partial _{k} , \overline{\partial}_{l} ).\\
\end{split}
\end{equation}
Since $\overline{\partial }_{S\otimes Q^{\ast}}\gamma =0$, we have  $D _{\overline{\partial}_{l}}^{H(t)}\gamma_{\bar{j}}=D _{\overline{\partial}_{j}}^{H(t)}\gamma_{\bar{l}}$, and
\begin{equation}
\begin{split}
 &g^{k\bar{l}}(D _{\partial _{k}}^{H(t)}D _{\overline{\partial}_{l}}^{H(t)}\gamma_{\bar{j}}+D _{\overline{\partial}_{l}}^{H(t)}D _{\partial _{k}}^{H(t)}\gamma_{\bar{j}})d\overline{z}^{j}\\
 =& 2g^{k\bar{l}}(D _{\partial _{k}}^{H(t)}D _{\overline{\partial}_{l}}^{H(t)}\gamma_{\bar{j}})d\overline{z}^{j}\\
 &- g^{k\bar{l}} (F_{H_{S}}(\partial _{k} , \overline{\partial}_{l} )\gamma_{\bar{j}} -\gamma_{\bar{j}}F_{H_{Q}}(\partial _{k} , \overline{\partial}_{l} ))d\overline{z}^{j}\\
 =& 2g^{k\bar{l}}(D _{\partial _{k}}^{H(t)}D _{\overline{\partial}_{j}}^{H(t)}\gamma_{\bar{l}})d\overline{z}^{j}\\
  & -(\sqrt{-1}\Lambda_{\omega} F_{H_{S}}\circ \gamma -\gamma \circ \sqrt{-1}\Lambda_{\omega}F_{H_{Q}})\\
 =& 2g^{k\bar{l}}(D _{\overline{\partial}_{j}}^{H(t)}D _{\partial _{k}}^{H(t)}\gamma_{\bar{l}})d\overline{z}^{j}\\
  & + 2g^{k\bar{l}} (F_{H_{S}}(\partial _{k} , \overline{\partial}_{j} )\gamma_{\bar{l}} -\gamma_{\bar{l}}F_{H_{Q}}(\partial _{k} , \overline{\partial}_{j} ))d\overline{z}^{j}\\
  & -(\sqrt{-1}\Lambda_{\omega}F_{H_{S}}\circ \gamma -\gamma \circ \sqrt{-1}\Lambda_{\omega}F_{H_{Q}}).\\
 \end{split}
\end{equation}
On the other hand, it is clear that
\begin{equation}
\begin{split}
& \overline{\partial }_{S\otimes Q^{\ast}}(\sqrt{-1}\Lambda_{\omega} \partial ^{H(t)}\gamma )=\overline{\partial }_{S\otimes Q^{\ast}}(g^{k\bar{l}}D _{\partial _{k}}^{H(t)}\gamma_{\bar{l}})\\
=& 2g^{k\bar{l}}(D _{\overline{\partial}_{j}}^{H(t)}D _{\partial _{k}}^{H(t)}\gamma_{\bar{l}})d\overline{z}^{j}+(D _{\partial _{k}}^{H(t)}\gamma_{\bar{l}})\frac{\partial g^{k\bar{l}}}{\partial \overline{z}^{j}} d\overline{z}^{j}\\
=& 2g^{k\bar{l}}(D _{\overline{\partial}_{j}}^{H(t)}D _{\partial _{k}}^{H(t)}\gamma_{\bar{l}})d\overline{z}^{j}+(D _{\partial _{k}}^{H(t)}\gamma_{\bar{j}}) g^{k\bar{l}}\nabla _{\overline{\partial}_{l}}d\overline{z}^{j}.\\
\end{split}
\end{equation}
The above equalities yield
\begin{equation}\label{r1}
\begin{split}
\Delta |\gamma (t)|_{H(t)}^{2}=& 2|\nabla ^{H(t)}\gamma |_{H(t)}^{2}+2Ric_{\omega}(\partial _{k}, \overline{\partial }_{j})g^{k\bar{l}}g^{i\bar{j}}\tr(\gamma_{\bar{l}}H_{Q}^{-1}\overline{(\gamma_{\bar{i}})^{T}}H_{S})\\
&+  4Re\{g^{k\bar{l}}\langle (F_{H_{S}}(\partial _{k} , \overline{\partial}_{j} )\gamma_{\bar{l}} -\gamma_{\bar{l}}F_{H_{Q}}(\partial _{k} , \overline{\partial}_{j} ))d\overline{z}^{j}, \gamma \rangle_{H(t)}\}\\
 & -2Re\{\langle(\sqrt{-1}\Lambda_{\omega}F_{H_{S}}\circ \gamma -\gamma \circ \sqrt{-1}\Lambda_{\omega}F_{H_{Q}}), \gamma \rangle_{H(t)}\}\\
 & +4Re\langle\overline{\partial }_{S\otimes Q^{\ast}}(\sqrt{-1}\Lambda_{\omega} \partial ^{H(t)}\gamma ), \gamma\rangle_{H(t)}.\\
\end{split}
\end{equation}
Combining (\ref{r0}), (\ref{S0}) and (\ref{Q0}), we have
\begin{equation}\label{r2}
\begin{split}
 \frac{\partial }{\partial t}|\gamma (t)|_{H(t)}^{2}=& \frac{\partial }{\partial t}g^{i\bar{j}}\tr(\gamma_{\bar{j}}H_{Q}^{-1}\overline{(\gamma_{\bar{i}})^{T}}H_{S})\\
 =& 4Re\langle\overline{\partial }_{S\otimes Q^{\ast}}(\sqrt{-1}\Lambda_{\omega} \partial ^{H(t)}\gamma ), \gamma\rangle_{H(t)}\\
 & -2Re\langle(\sqrt{-1}\Lambda_{\omega}(F_{H_{S}}-\gamma \wedge \gamma ^{\ast H(t)})-\lambda_{E}\Id_{S})\circ \gamma , \gamma \rangle_{H(t)}\\
 & +2Re\langle\gamma \circ (\sqrt{-1}\Lambda_{\omega}(F_{H_{Q}}-\gamma ^{\ast H(t)}\wedge \gamma )-\lambda_{E}\Id_{Q}), \gamma \rangle_{H(t)}.\\
\end{split}
\end{equation}
Then, from (\ref{r1}) and (\ref{r2}), we see that
\begin{equation}\label{r3}
\begin{split}
& (\Delta -\frac{\partial }{\partial t}) |\gamma (t)|_{H(t)}^{2}\\
=& 2|\nabla ^{H(t)}\gamma |_{H(t)}^{2}+2Ric_{\omega}(\partial _{k}, \overline{\partial }_{j})g^{k\bar{l}}g^{i\bar{j}}\tr(\gamma_{\bar{l}}H_{Q}^{-1}\overline{(\gamma_{\bar{i}})^{T}}H_{S})\\
&+4Re\{g^{k\bar{l}}\langle (F_{H_{S}}(\partial _{k} , \overline{\partial}_{j} )\gamma_{\bar{l}} -\gamma_{\bar{l}}F_{H_{Q}}(\partial _{k} , \overline{\partial}_{j} ))d\overline{z}^{j}, \gamma \rangle_{H(t)}\}\\
 & +2\langle(-\sqrt{-1}\Lambda_{\omega}\gamma \wedge \gamma ^{\ast H(t)})\circ \gamma +\gamma \circ (\sqrt{-1}\Lambda_{\omega}\gamma ^{\ast H(t)}\wedge \gamma), \gamma \rangle_{H(t)}.\\
 \end{split}
\end{equation}

By direct calculations, we obtain
\begin{equation}\label{G2}
\begin{split}
\Delta |G (t)|_{H(t)}^{2}=& 2g^{k\bar{l}}\partial_{k}\overline{\partial }_{l} |G (t)|_{H(t)}^{2}\\
=& 2g^{k\bar{l}}\{ \langle\nabla _{\partial_{k}}^{H(t)}\nabla _{\overline{\partial}_{l}}^{H(t)}G , G \rangle_{H(t)}+\langle G , \nabla _{\overline{\partial}_{k}}^{H(t)}\nabla _{\partial_{l}}^{H(t)}G \rangle_{H(t)}\\
& +\langle\nabla _{\partial_{k}}^{H(t)}G , \nabla _{\partial_{l}}^{H(t)}G \rangle_{H(t)} +\langle\nabla _{\overline{\partial}_{l}}^{H(t)}G , \nabla _{\overline{\partial}_{l}}^{H(t)}G \rangle_{H(t)}\}\\
=& 2Re \{g^{k\bar{l}} \langle\nabla _{\partial_{k}}^{H(t)}\nabla _{\overline{\partial}_{l}}^{H(t)}G +\nabla _{\overline{\partial}_{l}}^{H(t)}\nabla _{\partial_{k}}^{H(t)}G , G \rangle_{H(t)}\}\\
&+2|D_{H(t)}G |_{H(t)}^{2}\\
=& 4Re \{g^{k\bar{l}} \langle\nabla _{\partial_{k}}^{H(t)}\nabla _{\overline{\partial}_{l}}^{H(t)}G , G \rangle_{H(t)}\}+2|D_{H(t)}G |_{H(t)}^{2}\\
& -2Re\{\langle(\sqrt{-1}\Lambda_{\omega}F_{H_{S}}\circ G -G \circ \sqrt{-1}\Lambda_{\omega}F_{H_{Q}}), G \rangle_{H(t)}\},\\
\end{split}
\end{equation}
and
\begin{equation}\label{G3}
\begin{split}
\frac{\partial }{\partial t} |G (t)|_{H(t)}^{2}=& 2Re  \{\langle\frac{\partial }{\partial t}G , G \rangle_{H(t)}\}\\
& +Re\langle(H_{S}^{-1}\frac{\partial H_{S} }{\partial t})\circ G -G \circ (H_{Q}^{-1}\frac{\partial H_{Q} }{\partial t}), G \rangle_{H(t)}.\\
\end{split}
\end{equation}
Then we see (\ref{G0}), (\ref{S0}), (\ref{Q0}), (\ref{G2}) and (\ref{G3}) imply
\begin{equation}\label{G5}
\begin{split}
&(\Delta -\frac{\partial }{\partial t}) |G (t)|_{H(t)}^{2}\\
=& 4Re \{g^{k\bar{l}} \langle\nabla _{\partial_{k}}^{H(t)}(\nabla _{\overline{\partial}_{l}}^{H(t)}G-\gamma (\overline{\partial}_{l}))  , G \rangle_{H(t)}\}\\
&+2|D_{H(t)}G |_{H(t)}^{2}\\
& +2\langle(-\sqrt{-1}\Lambda_{\omega}\gamma \wedge \gamma ^{\ast H(t)})\circ G +G \circ (\sqrt{-1}\Lambda_{\omega}\gamma ^{\ast H(t)}\wedge \gamma), G \rangle_{H(t)},\\
\end{split}
\end{equation}

Using (\ref{G1}),  we have
\begin{equation}\label{G4}
\begin{split}
&(\Delta -\frac{\partial }{\partial t}) |G (t)|_{H(t)}^{2}\\
=& -4Re \{g^{k\bar{l}} \langle\nabla _{\partial_{k}}^{H(t)}\gamma_{0}(\overline{\partial}_{l})  , G \rangle_{H(t)}\}+2|D_{H(t)}G |_{H(t)}^{2}\\
& +2\langle(-\sqrt{-1}\Lambda_{\omega}\gamma \wedge \gamma ^{\ast H(t)})\circ G +G \circ (\sqrt{-1}\Lambda_{\omega}\gamma ^{\ast H(t)}\wedge \gamma), G \rangle_{H(t)}.\\
\end{split}
\end{equation}
Considering the second fundamental form $\gamma_{0}\in \Omega^{0,1}(S\otimes Q^{\ast})$, since $\overline{\partial }_{S\otimes Q^{\ast}}\gamma_{0} =0$, for any point $P\in M$, we have a domain $U_{P}$ and a local section  $G_{0}\in \Gamma (U_{P}; S\otimes Q^{\ast})$ such that $\gamma_{0}=\overline{\partial }_{S\otimes Q^{\ast}} G_{0}$. Locally, it holds that
\begin{equation}
\overline{\partial }_{S\otimes Q^{\ast}} (G+G_{0})=\gamma .
\end{equation}
Replacing $G$ by $G+G_{0}$ in (\ref{G5}), we see
\begin{equation}\label{G6}
\begin{split}
&(\Delta -\frac{\partial }{\partial t}) |G (t)+G_{0}|_{H(t)}^{2}\\
=& 4Re \{g^{k\bar{l}} \langle\nabla _{\partial_{k}}^{H(t)}(\nabla _{\overline{\partial}_{l}}^{H(t)}(G+G_{0})-\gamma (\overline{\partial}_{l}))  , G+G_{0} \rangle_{H(t)}\}\\
&+2|D_{H(t)}(G+G_{0}) |_{H(t)}^{2}\\
& +2\langle(-\sqrt{-1}\Lambda_{\omega}\gamma \wedge \gamma ^{\ast H(t)})\circ (G+G_{0}) , (G+G_{0}) \rangle_{H(t)}\\
&+2\langle(G+G_{0}) \circ (\sqrt{-1}\Lambda_{\omega}\gamma ^{\ast H(t)}\wedge \gamma), (G+G_{0}) \rangle_{H(t)}\\
=& 2|D_{H(t)}(G+G_{0}) |_{H(t)}^{2}\\
& +2\langle(-\sqrt{-1}\Lambda_{\omega}\gamma \wedge \gamma ^{\ast H(t)})\circ (G+G_{0}) , (G+G_{0}) \rangle_{H(t)}\\
&+2\langle(G+G_{0}) \circ (\sqrt{-1}\Lambda_{\omega}\gamma ^{\ast H(t)}\wedge \gamma), (G+G_{0}) \rangle_{H(t)}.\\
\end{split}
\end{equation}

\medskip

\section{$C^{0}$-estimate for the rescaled metrics}
\setcounter{equation}{0}

\medskip

Let $(M, \omega)$ be a compact K\"ahler manifold, and $\mathcal{E}=(E, \overline{\partial }_{E})$ be a nonsemistable holomorphic vector bundle on $M$.
We suppose that the length of the  HNS-filtration of $\mathcal{E}$ is
 one, i.e. there exists an exact sequence
\begin{equation}
0\rightarrow \mathcal{S}\rightarrow \mathcal{E} \rightarrow \mathcal{Q} \rightarrow 0,
\end{equation}
such that $\mathcal{S}$ is an $\omega$-stable  subsheaf, and $\mathcal{Q}$ is an  $\omega$-stable torsion-free coherent sheaf. We denote the singular set of $\mathcal{Q}$ by $\Sigma_{alg}$. Since $\mathcal{E}$ is nonsemistable, we have
\begin{equation}
\mu_{\omega }(\mathcal{S})>\mu_{\omega }(\mathcal{E})>\mu_{\omega }(\mathcal{Q}).
\end{equation}
Let $H(t)$ be the solution of the Hermitian-Yang-Mills flow (\ref{Dflow}) on the holomorphic bundle $(E, \overline{\partial }_{E})$ with initial metric $H_{0}$. As in the above section, we denote  $H_{\mathcal{S}}(t)$ and $H_{\mathcal{Q}}(t)$ the Hermitian metrics on $\mathcal{S}|_{M\setminus \Sigma_{alg}}$ and $\mathcal{Q}|_{M\setminus \Sigma_{alg}}$ induced by the metric $H(t)$ on $E$. In this section, we will derive a uniform local $C^{0}$-estimate for the rescaled metric $H_{\mathcal{S}}(t)=e^{2(\lambda_{\mathcal{S}}-\lambda_{\mathcal{E}})t}H_{\mathcal{S}}(t)$ and  $H_{\mathcal{Q}}(t)=e^{2(\lambda_{\mathcal{Q}}-\lambda_{\mathcal{E}})t}H_{\mathcal{Q}}(t)$ outside $\Sigma_{alg}$, where $\lambda_{\mathcal{S}}=\frac{\mu_{\omega }(\mathcal{S})}{\Vol(M, \omega)}$ and $\lambda_{\mathcal{Q}}=\frac{\mu_{\omega }(\mathcal{Q})}{\Vol(M, \omega)}$.

By Hironaka's flattening theorem  (\cite{Hi1}, \cite{Hi2}), we have a resolution of the  HNS-filtration (\cite{Sib}), by successively blowing up $\pi_{j}: M_{j}\rightarrow M_{j-1}$ with smooth center $Y_{i-1}\subset M_{i-1}$ finite times, where $j=1, \cdots , k$, $M_{0}=M$, $\tilde{M}=M_{k}$,  there is an exact sequence on $\tilde{M}$
\begin{equation}
0\rightarrow S\rightarrow \tilde{\mathcal{E}} \rightarrow Q \rightarrow 0
\end{equation}
such that: (1)$S, Q$ are locally free; (2)  the composition $\pi = \pi_{1}\circ \cdots \circ \pi_{k}: \tilde{M}\rightarrow M$ is biholomorphic outside $\Sigma_{alg} $; (3) $\tilde{\mathcal{E}}=\pi^{\ast}\mathcal{E}$; (4) $S, Q$ are isomorphic to the sheaves $\mathcal{S},  \mathcal{Q} $ respectively outside $\pi^{-1}\Sigma_{alg}$, $\pi_{\ast }S=\mathcal{S}$ and $\mathcal{Q}^{\ast \ast }=(\pi_{\ast }Q )^{\ast \ast}$.

It is well known that $\tilde{M}$ is also K\"ahler (\cite{GH}). As in \cite{BS}, we fix arbitrary K\"ahler metrics $\eta_{i} $ on $M_{i}$ and set
 \begin{equation}\label{omegai}
 \omega_{1, \epsilon }=\pi_{1}^{\ast }\omega +\epsilon_{1}\eta_{1},  \quad  \quad  \omega_{i, \epsilon }=\pi_{i}^{\ast }\omega_{i-1 , \epsilon } +\epsilon_{i}\eta_{i}
 \end{equation}
for all $1\leq i \leq k$ and $0<\epsilon_{i} \leq 1$, where $\epsilon =(\epsilon_{1}, \cdots , \epsilon_{k})$.
For the sake of simplicity, we only consider the case that $k=1$,  the general case following by induction.

In the following, we always assume that $\pi : \tilde{M}\rightarrow M$ is a single blow-up with smooth centre, fix a K\"ahler metric $\eta $ on $\tilde{M}$ and set $\omega_{\epsilon }=\pi^{\ast}\omega +\epsilon \eta $
  for $0<\epsilon \leq 1$. In fact there exists a holomorphic line bundle $L$ over $\tilde{M}$ with respect to the exceptional divisor $D\subset \pi^{-1}(\Sigma_{alg})$ such that the $(1, 1)$-form $\pi^{\ast }\omega +\delta \sqrt{-1}F_{H_{L}}$ is positive for some $\delta$ small enough (for the proof see for example \cite{GH,VO}), where $\sqrt{-1}F_{H_{L}}$ is the Chern form with respect to some Hermitian metric on $L$. In the following, we can set the K\"ahler metric $\eta $ by
   \begin{equation}\label{Def05}
   \eta =\pi^{\ast }\omega +\delta \cdot \sqrt{-1}F_{H_{L}}.
   \end{equation}
    Bando and Siu (Lemma 3 in \cite{BS}) derived a uniform  Sobolev inequality for $(\tilde{M}, \omega_{\epsilon})$, i.e. there exists a uniform constant $C_{S}$ such that
  \begin{equation}\label{Sobolev}
(\int_{\tilde{M}} |\rho|^{\frac{2n}{2n-1}}\frac{\omega_{\epsilon}^{n}}{n!})^{\frac{2n-1}{2n}}\leq C_{S}\Big( \int_{\tilde{M}} |d\rho |_{\omega_{\epsilon }} +|\rho|\frac{\omega_{\epsilon}^{n}}{n!} \Big)
\end{equation}
for all $\rho \in C^{1}(\tilde{M})$ and all $0<\epsilon \leq 1$. Using Li's result (Proposition 3 in \cite{Pli}), we obtain a uniform lower bound on the first eigenvalue of $\Delta_{\epsilon }$,  i.e. there exists a uniform constant $C_{P}$ such that
  \begin{equation}\label{Poincare}
C_{P}\cdot \inf_{a\in R}(\int_{\tilde{M}} |\rho-a|^{2}\frac{\omega_{\epsilon}^{n}}{n!})\leq  \int_{\tilde{M}} |d\rho |_{\omega_{\epsilon }}^{2} \frac{\omega_{\epsilon}^{n}}{n!}
\end{equation}
for all $\rho \in W^{1, 2}(\tilde{M})$ and all $0<\epsilon \leq 1$.
 Combining Cheng and Li's estimate (\cite{CL}) with  Grigor'yan's result (Theorem 1.1 in \cite{Gr}), we have the following uniform upper bounds of the heat kernels and  uniform lower bounds of the Green functions.

 \medskip

\begin{proposition}\label{Prop 2.1}{\bf(Proposition 2 in \cite{BS})}
Let $\mathcal{K}_{\epsilon}$ be the heat kernel with respect to the metric $\omega_{\epsilon}$, then for any $\tau >0$, there exists a constant $C_{\mathcal{K}}(\tau)$ independent of $\epsilon $, such that
\begin{equation}\label{kernel01}0\leq \mathcal{K}_{\epsilon}(x , y, t)\leq C_{\mathcal{K}}(\tau) (t^{-n}\exp{(-\frac{(d_{\omega_{\epsilon}}(x, y))^{2}}{(4+\tau )t})}+1)\end{equation} for every $x, y\in \tilde{M}$ and $0<t < +\infty$, where $d_{\omega_{\epsilon}}(x, y)$ is the distance between $x$ and $y$ with respect to the metric $\omega_{\epsilon}$. There also exists a constant $C_{\mathcal{G}}$ such that
\begin{equation}\label{green} \mathcal{G}_{\epsilon}(x , y)\geq -C_{\mathcal{G}}\end{equation}
for every $x, y\in \tilde{M}$ and $0<\epsilon \leq 1$, where $\mathcal{G}_{\epsilon}$ is the Green function with respect to the metric $\omega_{\epsilon}$.
\end{proposition}

 \medskip

Let $H_{\epsilon}(t)$ be the solution of the Hermitian-Yang-Mills flow (\ref{Dflow}) on the holomorphic bundle $\tilde{\mathcal{E}}$ over the K\"ahler manifold $(\tilde{M}, \omega_{\epsilon})$ with the fixed initial metric $\pi^{\ast}H_{0}$, i.e. it satisfies
\begin{equation}\label{SSS1}
\left \{\begin{split} &H_{\epsilon}(t)^{-1}\frac{\partial H_{\epsilon}(t)}{\partial
t}=-2(\sqrt{-1}\Lambda_{\omega_{\epsilon}}F_{H_{\epsilon}(t)}-\lambda_{\mathcal{E}, \epsilon } \Id_{\tilde{\mathcal{E}}}),\\
&H_{\epsilon}(0)=\pi^{\ast}H_{0}.\\
\end{split}
\right.
\end{equation}
By Bando and Siu's argument in \cite{BS} and the uniqueness of the Hermitian-Yang-Mills flow, we know that, by choosing a subsequence, $H_{ \epsilon}(t)$ converges to $H(t)$ in $C_{loc}^{\infty}$-topology outside $\Sigma_{alg }$ as $\epsilon \rightarrow 0$. We also denote  $H_{S, \epsilon}(t)$ and $H_{Q, \epsilon }(t)$ the Hermitian metrics on bundles $S$ and $Q$ induced by the metric $H_{\epsilon}(t)$. It is easy to see that: by choosing a subsequence,
\begin{equation}
H_{S, \epsilon}(t)\rightarrow H_{\mathcal{S}}(t), \quad H_{Q, \epsilon }(t)\rightarrow H_{\mathcal{Q} }(t)
\end{equation}
in $C_{loc}^{\infty}$-topology outside $\Sigma_{alg }$ as $\epsilon \rightarrow 0$.

Since $\mathcal{S}$ and $\mathcal{Q}$ are $\omega$-stable, it is easy to check that $(S, \overline{\partial }_{S})$ and $(Q, \overline{\partial }_{Q})$ are $\omega_{\epsilon }$-stable for sufficiently small $\epsilon$. By Donaldson-Uhlenbeck-Yau theorem,   we can suppose that  $K_{S, \epsilon }$ and $K_{Q, \epsilon }$ are $\omega_{\epsilon}$-Hermitian-Einstein metrics on  $(S, \overline{\partial }_{S})$ and $(Q, \overline{\partial }_{Q})$, i.e.
\begin{equation}
\sqrt{-1}\Lambda_{\omega_{\epsilon} }F_{K_{S, \epsilon}}=\lambda_{S, \epsilon}\Id_{S}
\end{equation}
and
\begin{equation}
\sqrt{-1}\Lambda_{\omega_{\epsilon} }F_{K_{Q, \epsilon}}=\lambda_{Q, \epsilon}\Id_{Q}.
\end{equation}
Here $\lambda_{S, \epsilon}=\frac{\mu_{\omega_{\epsilon} }(S)}{\Vol(\tilde{M}, \omega_{\epsilon})}$ and $\lambda_{Q, \epsilon}=\frac{\mu_{\omega_{\epsilon} }(Q)}{\Vol(\tilde{M}, \omega_{\epsilon})}$, and
\begin{equation}
\lambda_{S, \epsilon}\rightarrow \lambda_{\mathcal{S}}, \quad \lambda_{Q, \epsilon}\rightarrow \lambda_{\mathcal{Q}}
\end{equation}
as $\epsilon \rightarrow 0$.

 Denote $h_{S, \epsilon}(t)=K_{S, \epsilon}^{-1}H_{S, \epsilon}(t)$, $h_{Q, \epsilon}(t)=K_{Q, \epsilon}^{-1}H_{Q, \epsilon}(t)$, and set $\tilde{h}_{S, \epsilon}(t)=e^{2(\lambda_{S, \epsilon}-\lambda_{\tilde{\mathcal{E}}, \epsilon })t}h_{S}(t)$, $\tilde{h}_{Q, \epsilon}(t)=e^{2(\lambda_{Q, \epsilon}-\lambda_{\tilde{\mathcal{E}}, \epsilon })t}h_{Q, \epsilon}(t)$. Using (\ref{S0}) and (\ref{Q0}), we have
\begin{equation}\label{TTTT1}
\begin{split}
(\Delta_{\epsilon} -\frac{\partial }{\partial t})\tr \tilde{h}_{S, \epsilon}=&2\tr (-\sqrt{-1}\Lambda _{\omega_{\epsilon} } \bar{\partial }\tilde{h}_{S, \epsilon} \tilde{h}_{S, \epsilon}^{-1}\partial _{K_{S, \epsilon}}\tilde{h}_{S, \epsilon})-2(\lambda_{S, \epsilon}-\lambda_{\tilde{\mathcal{E}}, \epsilon})\tr(\tilde{h}_{S, \epsilon})\\
&-2\tr (\tilde{h}_{S, \epsilon}\sqrt{-1}\Lambda_{\omega_{\epsilon} }(F_{H_{S, \epsilon}}-F_{K_{S, \epsilon}}))
-\tr(\tilde{h}_{S, \epsilon} H_{S, \epsilon}^{-1}\frac{\partial H_{S, \epsilon}}{\partial t})\\
=&2\tr (-\sqrt{-1}\Lambda _{\omega_{\epsilon} } \bar{\partial }\tilde{h}_{S, \epsilon} \tilde{h}_{S, \epsilon}^{-1}\partial _{K_{S, \epsilon}}\tilde{h}_{S, \epsilon})\\
&-2\tr (\tilde{h}_{S, \epsilon}(\sqrt{-1}\Lambda_{\omega_{\epsilon} }\gamma_{\epsilon} \wedge \gamma_{\epsilon}^{\ast }))\\
\geq  & 0,\\
\end{split}
\end{equation}
and
\begin{equation}\label{QQQQ1}
\begin{split}
(\Delta_{\epsilon} -\frac{\partial }{\partial t})\tr \tilde{h}_{Q, \epsilon}^{-1}=&2\tr (-\sqrt{-1}\Lambda _{\omega_{\epsilon} } \bar{\partial }\tilde{h}_{Q, \epsilon}^{-1} \circ \tilde{h}_{Q, \epsilon}\circ  \partial _{H_{Q, \epsilon}}\tilde{h}_{Q, \epsilon}^{-1})\\
&+2\tr (\tilde{h}_{Q, \epsilon}^{-1}(\sqrt{-1}\Lambda_{\omega_{\epsilon} }\gamma_{\epsilon}^{\ast }\wedge \gamma_{\epsilon}  ))\\
\geq  & 0,\\
\end{split}
\end{equation}
where we have used the nonnegativity of $\sqrt{-1}\Lambda_{\omega_{\epsilon} }\gamma_{\epsilon}^{\ast }\wedge \gamma_{\epsilon}$.
Using the above inequalities and the uniform upper bound on the heat kernels, we can get a uniform bound (independent of $\epsilon $) on $\tr \tilde{h}_{S, \epsilon } + \tr \tilde{h}_{Q, \epsilon }^{-1}$. In fact, we obtain the following lemma.

\medskip

\begin{lemma}\label{l4.2}
There exists a uniform constant $C_{0, 1}$ such that
\begin{equation}\label{C00}
\sup_{x\in \tilde{M}}(\tr \tilde{h}_{S, \epsilon }(x, t) + \tr \tilde{h}_{Q, \epsilon }^{-1}(x, t))\leq C_{0, 1},
\end{equation}
for all $t\in [0, +\infty )$ and $0<\epsilon \leq 1$.
\end{lemma}

\medskip

{\bf Proof. } In \cite{BS}, Bando and Siu conjectured that $K_{S, \epsilon}$ and $K_{Q, \epsilon}$ should converge to $\omega$-Hermitian-Einstein metrics $K_{\mathcal{S}}$ and $K_{\mathcal{S}}$ on sheaves $\mathcal{S}$ and $\mathcal{Q}$ in local $C^{\infty}$-topology outside  $\pi^{-1}\Sigma_{alg} $. This fact has been proved in \cite{LZZ}. By taking  constants on $K_{S, \epsilon}$ and $K_{Q, \epsilon}$, we can suppose that
 \begin{equation}\label{det2}
 \int_{\tilde{M}}\log \det(\tilde{h}_{S, \epsilon}(0))\frac{\omega_{\epsilon}^{n}}{n!}= \int_{\tilde{M}}\log \det(\tilde{h}_{Q, \epsilon}(0))\frac{\omega_{\epsilon}^{n}}{n!}=0 .
 \end{equation}
By the uniform $L^{1}$-estimate in \cite{LZZ} (Lemma 5.1.), we see that there exists a uniform constant $\hat{C}$ such that
 \begin{equation}\label{det20}
 \int_{\tilde{M}}\log \{\tr(\tilde{h}_{S, \epsilon}(0))+\tr(\tilde{h}_{S, \epsilon}^{-1}(0))\}+\log \{\tr(\tilde{h}_{Q, \epsilon}(0))+\tr(\tilde{h}_{Q, \epsilon}^{-1}(0))\}\frac{\omega_{\epsilon}^{n}}{n!}\leq \hat{C}
 \end{equation}
 for all $0<\epsilon \leq 1$.
For any point $x\in \tilde{M}\setminus \pi^{-1}(\Sigma_{alg})$ and any vector $X\in \tilde{\mathcal{E}}_{x}$ satisfying  $|X|_{\pi^{\ast }H_{0}}=1$, it is easy to check that
\begin{equation}
-\sup_{x\in M}|F_{H_{0}}|_{\omega , H_{0}}(x) \cdot \pi^{\ast }\omega \leq \langle   \sqrt{-1} F_{\pi^{\ast }H_{0}}(X), X\rangle_{\pi^{\ast }H_{0}} \leq \sup_{x\in M}|F_{H_{0}}|_{\omega , H_{0}}(x) \cdot \pi^{\ast }\omega ,
\end{equation}
and
\begin{equation}
-n\sup_{x\in M}|F_{H_{0}}|_{\omega , H_{0}}(x)  \leq \Lambda_{\omega_{\epsilon}}\langle   \sqrt{-1} F_{\pi^{\ast }H_{0}}(X), X\rangle_{\pi^{\ast }H_{0}} \leq n\sup_{x\in M}|F_{H_{0}}|_{\omega , H_{0}}(x) .
\end{equation}
Then we have
\begin{equation}\label{bbin}
-n\sup_{x\in M}|F_{H_{0}}|_{\omega , H_{0}}(x) \Id_{\tilde{\mathcal{E}}} \leq \sqrt{-1}\Lambda_{\omega_{\epsilon}}  F_{\pi^{\ast }H_{0}} \leq n\sup_{x\in M}|F_{H_{0}}|_{\omega , H_{0}}(x) \Id_{\tilde{\mathcal{E}}}
\end{equation}
on the whole $\tilde{M}$. Using the Gauss-Codazzi equation (\ref{GC}), we have
\begin{equation}\label{FS}
 \sqrt{-1}\Lambda_{\omega_{\epsilon}}  F_{H_{S, \epsilon}(0)} \leq n\sup_{x\in M}|F_{H_{0 }}|_{\omega , H_{0}}(x) \Id_{S},
\end{equation}
and
\begin{equation}\label{FQ}
 \sqrt{-1}\Lambda_{\omega_{\epsilon}}  F_{H_{Q, \epsilon}(0)} \geq -n\sup_{x\in M}|F_{H_{0 }}|_{\omega , H_{0}}(x) \Id_{Q}.
\end{equation}
Direct calculations yield
\begin{equation}\label{TTTT10}
\begin{split}
\Delta_{\epsilon} \tr \tilde{h}_{S, \epsilon}(0)=&2\tr (-\sqrt{-1}\Lambda _{\omega_{\epsilon} } \bar{\partial }\tilde{h}_{S, \epsilon}(0) \tilde{h}_{S, \epsilon}^{-1}(0)\partial _{K_{S, \epsilon}}\tilde{h}_{S, \epsilon}(0))\\
&-2\tr (\tilde{h}_{S, \epsilon}(0)\sqrt{-1}\Lambda_{\omega_{\epsilon} }(F_{H_{S, \epsilon}(0)}-F_{K_{S, \epsilon}}))
\\
\geq &2\tr (-\sqrt{-1}\Lambda _{\omega_{\epsilon} } \bar{\partial }\tilde{h}_{S, \epsilon}(0) \tilde{h}_{S, \epsilon}^{-1}(0)\partial _{K_{S, \epsilon}}\tilde{h}_{S, \epsilon}(0))\\
&-2\tr (\tilde{h}_{S, \epsilon}(0))(n\sup_{x\in M}|F_{H_{0 }}|_{\omega , H_{0}}(x) -\lambda_{S, \epsilon}),
\\
\end{split}
\end{equation}
\begin{equation}\label{TTTT101}
\begin{split}
\Delta_{\epsilon} \tr \tilde{h}_{Q, \epsilon}^{-1}(0)\geq &2\tr (-\sqrt{-1}\Lambda _{\omega_{\epsilon} } \bar{\partial }\tilde{h}_{Q, \epsilon}^{-1} \circ \tilde{h}_{Q, \epsilon}\circ  \partial _{H_{Q, \epsilon}}\tilde{h}_{Q, \epsilon}^{-1})\\
&-2\tr (\tilde{h}_{Q, \epsilon}^{-1}(0))(n\sup_{x\in M}|F_{H_{0 }}|_{\omega , H_{0}}(x) +\lambda_{Q, \epsilon}).
\\
\end{split}
\end{equation}
Then there is a uniform constant $\hat{C}_{0, 1}$ such that
\begin{equation}\label{TTTT1001}
\Delta_{\epsilon} \log (\tr \tilde{h}_{S, \epsilon}(0)+ 1) \geq   -\hat{C}_{0, 1},
\end{equation}
and
\begin{equation}\label{TTTT1002}
\Delta_{\epsilon} \log (\tr \tilde{h}_{Q, \epsilon}^{-1}(0)+ 1) \geq   -\hat{C}_{0, 1},
\end{equation}
for all $0<\epsilon \leq 1$. The uniform lower  bounds of the Green functions (\ref{green}), inequalities (\ref{det20}), (\ref{TTTT1001}) and (\ref{TTTT1002}) imply that there is a uniform constant $\tilde{C}_{0, 1}$ such that
\begin{equation}\label{TTTT1003}
\sup_{x\in \tilde{M}} (\tr \tilde{h}_{S, \epsilon}(0)+ \tr \tilde{h}_{Q, \epsilon}^{-1}(0)) \leq   \tilde{C}_{0, 1},
\end{equation}
for all $0<\epsilon \leq 1$. By (\ref{TTTT1}), (\ref{QQQQ1}) and the maximum principle, we obtain the estimate (\ref{C00}).

\hfill $\Box$ \\

In the following we will get uniform upper bounds on $\det(\tilde{h}_{S, \epsilon }^{-1}(t))$ and $\det(\tilde{h}_{Q, \epsilon }(t))$, which imply uniform upper bounds on $\tr \tilde{h}_{S, \epsilon }^{-1}(t)$ and $ \tr \tilde{h}_{Q, \epsilon }(t)$ by (\ref{C00}). Then we get uniform $C^{0}$-bounds on $\tilde{h}_{S, \epsilon }(t)$ and  $\tilde{h}_{Q, \epsilon }(t)$.

Let's recall
\begin{equation}
(\Delta_{\epsilon} -\frac{\partial }{\partial t})\tr (\sqrt{-1}\Lambda_{\omega_{\epsilon } }F_{H_{\epsilon} (t)} -\lambda_{\mathcal{E}, \epsilon } \Id_{\mathcal{E}}) =0.
\end{equation}
For simplicity, we set $f_{\epsilon}(t)=\tr (\sqrt{-1}\Lambda_{\omega_{\epsilon } }F_{H_{\epsilon} (t)} -\lambda_{\mathcal{E}, \epsilon } \Id_{\mathcal{E}})$. By (\ref{bbin}) and the uniform Poincar\'e inequality (\ref{Poincare}), we know that there exists a uniform constant $C_{F, 1}$ such that
\begin{equation}
\int_{\tilde{M}} |f_{\epsilon }(t)|^{2}\frac{\omega_{\epsilon}^{n}}{n!}\leq C_{F, 1}\cdot \exp\{-2C_{p}t\}
\end{equation}
for all $0<\epsilon \leq 1$ and $t\in [0, +\infty )$. From the upper bound of the heat kernels (\ref{kernel01}),  we get
\begin{equation}\label{FF1}
 \max_{x\in \tilde{M}}f_{\epsilon}^{2}(x, t) \leq 2C_{\mathcal{K}}(\tau )\int_{\tilde{M}}|f_{\epsilon }(t-1)|^{2}\frac{\omega_{\epsilon }^{n}}{n!}\leq 2C_{\mathcal{K}}(\tau ) \cdot C_{F, 1}\cdot \exp\{-2C_{p}t\}
\end{equation}
for all  all $0<\epsilon \leq 1$ and $t\in [1, +\infty )$.
Using the Gauss-Codazzi equation (\ref{GC}), we have
\begin{equation}\label{C04}
\begin{split}
& \tr (\sqrt{-1}\Lambda_{\omega_{\epsilon } }F_{H_{\epsilon} (t)} -\lambda_{\mathcal{E}, \epsilon } \Id_{\mathcal{E}})\\
=& \tr (\sqrt{-1}\Lambda_{\omega_{\epsilon } }(F_{H_{S, \epsilon } (t)}-\gamma_{\epsilon} \wedge \gamma_{\epsilon}^{\ast})) -\lambda_{S, \epsilon}\rank (S)\\
& +\tr (\sqrt{-1}\Lambda_{\omega_{\epsilon} }(F_{H_{Q, \epsilon} (t)})-\gamma_{\epsilon}^{\ast} \wedge \gamma_{\epsilon})-\lambda_{Q, \epsilon } \rank (Q),\\
\end{split}
\end{equation}
and then
\begin{equation}\label{C01}
\begin{split}
& |\sqrt{-1}\Lambda_{\omega_{\epsilon } }(F_{H_{S, \epsilon } (t)}-\gamma_{\epsilon} \wedge \gamma_{\epsilon}^{\ast}))-\lambda_{\mathcal{E} , \epsilon}\Id_{S}|^{2}\\
&+|\sqrt{-1}\Lambda_{\omega_{\epsilon} }(F_{H_{Q, \epsilon} (t)})-\gamma_{\epsilon}^{\ast} \wedge \gamma_{\epsilon})-\lambda_{\mathcal{E} , \epsilon}\Id_{Q}|^{2}\\
=& 2(\lambda_{S, \epsilon }-\lambda_{Q, \epsilon}) \tr (\sqrt{-1}\Lambda_{\omega_{\epsilon } }(F_{H_{S, \epsilon } (t)}-\gamma_{\epsilon} \wedge \gamma_{\epsilon}^{\ast})) -\lambda_{S, \epsilon}\Id_{S})\\
& + 2(\lambda_{Q, \epsilon }-\lambda_{\mathcal{E}, \epsilon }) \tr (\sqrt{-1}\Lambda_{\omega_{\epsilon } }F_{H_{\epsilon} (t)} -\lambda_{\mathcal{E}, \epsilon } \Id_{\mathcal{E}})\\
& +|\sqrt{-1}\Lambda_{\omega_{\epsilon } }(F_{H_{S, \epsilon } (t)}-\gamma_{\epsilon} \wedge \gamma_{\epsilon}^{\ast})) -\lambda_{S, \epsilon}\Id_{S}|^{2}\\
& +|\sqrt{-1}\Lambda_{\omega_{\epsilon} }(F_{H_{Q, \epsilon} (t)})-\gamma_{\epsilon}^{\ast} \wedge \gamma_{\epsilon})-\lambda_{Q, \epsilon }\Id_{Q}|^{2}\\
& +(\lambda_{S, \epsilon }-\lambda_{\mathcal{E}, \epsilon })^{2}\rank(S) + (\lambda_{Q, \epsilon }-\lambda_{\mathcal{E}, \epsilon })^{2}\rank(Q).\\
\end{split}
\end{equation}

\medskip

Let's recall the Donaldson's functional
 \begin{equation}
 \mathcal{M} _{\tilde{\mathcal{E}}, \epsilon}^{0}(H_{0}, H)=\int_{0}^{1}\int_{\tilde{M}} \tr(\sqrt{-1}\Lambda_{\omega }F_{H(s)}H^{-1}(s)\frac{\partial H(s)}{\partial s})\frac{\omega_{\epsilon}^{n}}{n!}
 \end{equation}
 and
 \begin{equation}
  \mathcal{M} _{\tilde{\mathcal{E}}, \epsilon }(H_{0}, H)= \mathcal{M} _{\tilde{\mathcal{E}}}^{0}(H_{0}, H)-\lambda_{\tilde{\mathcal{E}}, \epsilon }\int_{\tilde{M}} \log \det (H_{0}^{-1}H)\frac{\omega_{\epsilon}^{n}}{n!},
 \end{equation}
 where $H(s)$ is a path connecting metrics $H_{0}$ and $H$ on $\tilde{\mathcal{E}}$. Donaldson proved that the above integral is independent of the path, and the Hermitian-Yang-Mills flow is the gradient flow of the above functional, i.e. if $H_{\epsilon}(t)$ is a solution of the Hermitian-Yang-Mills flow (\ref{Dflow}),  we have
 \begin{equation}\label{M0}
 \frac{d}{dt} \mathcal{M} _{\tilde{\mathcal{E}}, \epsilon}(H_{0}, H(t))=-2\int_{\tilde{M}} |\sqrt{-1}\Lambda_{\omega_{\epsilon } }F_{H_{\epsilon} (t)} -\lambda_{\tilde{\mathcal{E}}, \epsilon } \Id_{\mathcal{E}}|^{2}_{H_{\epsilon}(t)}\frac{\omega_{\epsilon}^{n}}{n!}.
 \end{equation}

\medskip

Since $(S, \overline{\partial }_{S})$ and $(Q, \overline{\partial }_{Q})$ are $\omega_{\epsilon }$-stable bundles over $\tilde{M}$ for a fixed $\epsilon$, the Donaldson's functional $\mathcal{M} _{S, \epsilon} (H_{S, 0}, \cdot )$ and $\mathcal{M} _{Q, \epsilon} (H_{Q, 0}, \cdot )$ are bounded from below. In the following, we will prove that the Donaldson's functional $\mathcal{M} _{S, \epsilon} (H_{S, \epsilon }(0), \cdot )$ and $\mathcal{M} _{Q, \epsilon} (H_{Q, \epsilon }(0), \cdot )$ are uniformly bounded from below.

\medskip

\begin{proposition}
 There exists a uniform positive constant $C_{\mathcal{M}}$ such that
 \begin{equation}\label{MS}
 \mathcal{M} _{S, \epsilon} (H_{S, \epsilon }(0), K_{S, \epsilon} )\geq -C_{\mathcal{M}}
 \end{equation}
 and
 \begin{equation}\label{MQ}
 \mathcal{M} _{Q, \epsilon} (H_{Q, \epsilon }(0), K_{Q, \epsilon} )\geq -C_{\mathcal{M}}
 \end{equation}
for all $0<\epsilon \leq 1$, where $H_{S, \epsilon }(0)$ and $H_{Q, \epsilon }(0)$ are metrics on $S$ and $Q$ induced by the pull back metric $\pi^{\ast }H_{0}$, $K_{S, \epsilon }$ and $K_{Q, \epsilon }$ are $\omega_{\epsilon}$-Hermitian-Einstein metrics on $S$ and $Q$.
\end{proposition}

\medskip

{\bf Proof. } Setting $\exp \{\varrho_{S, \epsilon }\}:=h_{S, \epsilon}^{-1}(0)=H_{S, \epsilon}^{-1}(0)K_{S, \epsilon}$, $\exp \{\varrho_{Q, \epsilon }\}:=h_{Q, \epsilon}^{-1}(0)=H_{Q, \epsilon}^{-1}(0)K_{Q, \epsilon}$, we have the following expression of the Donaldson's functional
\begin{equation}\label{7}
\begin{split}
\mathcal{M} _{S, \epsilon} (H_{S, \epsilon }(0), K_{S, \epsilon} )  = & \int_{\tilde{M} } \tr \{\varrho_{S, \epsilon } \cdot(\sqrt{-1}\Lambda_{\omega_{\epsilon} }F_{H_{S, \epsilon }(0) }-\lambda_{S, \epsilon }\Id_{S})\}\\ &+ \langle \Psi (\varrho_{S, \epsilon }) (\overline{\partial } \varrho_{S, \epsilon }) , \overline{\partial} \varrho_{S, \epsilon } \rangle_{H_{S, \epsilon }(0)}\frac{\omega_{\epsilon}^{n}}{n!},
\end{split}
\end{equation}
where $\Psi (x, y)= (x-y)^{-2}(e^{y-x }-(y-x)-1)$. Since the second part of the right hand side of the above equality is nonnegative, we have
\begin{equation}\label{70}
\mathcal{M} _{S, \epsilon} (H_{S, \epsilon }(0), K_{S, \epsilon} ) \geq  \int_{\tilde{M} } \tr \{\varrho_{S, \epsilon } \cdot(\sqrt{-1}\Lambda_{\omega_{\epsilon} }F_{H_{S, \epsilon }(0) }-\lambda_{S, \epsilon }\Id_{S})\}\frac{\omega_{\epsilon}^{n}}{n!},
\end{equation}
and
\begin{equation}\label{q70}
\mathcal{M} _{Q, \epsilon} (H_{Q, \epsilon }(0), K_{Q, \epsilon} ) \geq  \int_{\tilde{M} } \tr \{\varrho_{Q, \epsilon } \cdot(\sqrt{-1}\Lambda_{\omega_{\epsilon} }F_{H_{Q, \epsilon }(0) }-\lambda_{Q, \epsilon }\Id_{S})\}\frac{\omega_{\epsilon}^{n}}{n!}.
\end{equation}
By (\ref{TTTT1003}), there must exist a uniform constant $ \tilde{C}_{0, 2}$ such that
\begin{equation}\label{SQ}
\varrho_{S, \epsilon } \geq -\tilde{C}_{0, 2}\Id_{S}, \quad \varrho_{Q, \epsilon } \leq \tilde{C}_{0, 2}\Id_{S}
\end{equation}
for all $0<\epsilon \leq 1$. (\ref{det20}), (\ref{FS}), (\ref{SQ}) and (\ref{70}) imply that there is a uniform $ \tilde{C}_{0, 3}$ such that
\begin{equation}\label{700}
\mathcal{M} _{S, \epsilon} (H_{S, \epsilon }(0), K_{S, \epsilon} ) \geq  \int_{\tilde{M} } \tr (\varrho_{S, \epsilon }) \cdot \tr(\sqrt{-1}\Lambda_{\omega_{\epsilon} }F_{H_{S, \epsilon }(0) })\frac{\omega_{\epsilon}^{n}}{n!}-\tilde{C}_{0, 3},
\end{equation}
and
\begin{equation}\label{q700}
\mathcal{M} _{Q, \epsilon} (H_{Q, \epsilon }(0), K_{Q, \epsilon} ) \geq  \int_{\tilde{M} } \tr (\varrho_{Q, \epsilon }) \cdot \tr(\sqrt{-1}\Lambda_{\omega_{\epsilon} }F_{H_{Q, \epsilon }(0) })\frac{\omega_{\epsilon}^{n}}{n!}-\tilde{C}_{0, 3},
\end{equation}
for all $0<\epsilon \leq 1$.

 By the definition, we have
\begin{equation}\label{SQ1}
\tr \varrho_{S, \epsilon } =\log \frac{\det (K_{S, \epsilon })}{ \det (H_{S, \epsilon }(0))}.
\end{equation}
In the following, we will show that $\epsilon \cdot \log \frac{\det (K_{S, \epsilon })}{ \det (H_{S, \epsilon }(0))}$ are uniformly bounded. Let $\tilde{H}_{S, \epsilon }(t )$ be the evolved metric along the Hermitian-Yang-Mills flow (\ref{Dflow}) with initial metric $H_{S, \epsilon }(0)$, i.e. it satisfies
\begin{equation}\label{SSS11}
\left \{\begin{split} &\tilde{H}_{S, \epsilon}(t)^{-1}\frac{\partial \tilde{H}_{S, \epsilon}(t)}{\partial
t}=-2(\sqrt{-1}\Lambda_{\omega_{\epsilon}}F_{\tilde{H}_{S, \epsilon}(t)}-\lambda_{S, \epsilon } \Id_{S}),\\
&\tilde{H}_{S, \epsilon}(0)=H_{S, \epsilon }(0).\\
\end{split}
\right.
\end{equation}
Since $H_{S, \epsilon }(0)$ is a fixed smooth Hermitian metric on bundle $S$, there is a uniform constant $\tilde{C}_{0, 5}$ such that
\begin{equation}\label{80}
-\tilde{C}_{0, 5}\eta \leq \sqrt{-1}\tr F_{H_{S, \epsilon }(0)}\leq \tilde{C}_{0, 5} \eta ,
\end{equation}
for all $0<\epsilon \leq 1$. (\ref{80}) and (\ref{SSS11}) imply that
\begin{equation}
\begin{split}
|\epsilon \tr (\sqrt{-1}\Lambda_{\omega_{\epsilon}}F_{H_{S, \epsilon }(0)})|&=n\epsilon |\frac{\sqrt{-1}\tr F_{H_{S, \epsilon }(0)} \wedge \omega_{\epsilon}^{n-1}}{\omega_{\epsilon}^{n}}|\\
&\leq n\frac{\epsilon \eta \wedge \omega_{\epsilon}^{n-1}}{\omega_{\epsilon}^{n}}\leq n \tilde{C}_{0, 5},\\
\end{split}
\end{equation}
and
\begin{equation}\label{dd01}
\begin{split}
\epsilon |\frac{\partial }{\partial t}\log \frac{\det (\tilde{H}_{S, \epsilon }(t))}{ \det (H_{S, \epsilon }(0))}|&=\epsilon |\tr (\tilde{H}_{S, \epsilon}(t)^{-1}\frac{\partial \tilde{H}_{S, \epsilon}(t)}{\partial
t})|\\ &=|\epsilon \tr (\sqrt{-1}\Lambda_{\omega_{\epsilon}}F_{H_{S, \epsilon }(t)}-\lambda_{S, \epsilon } \Id_{S})|\\&\leq
\sup_{x\in \tilde{M}}|\epsilon \tr (\sqrt{-1}\Lambda_{\omega_{\epsilon}}F_{H_{S, \epsilon }(0)}-\lambda_{S, \epsilon } \Id_{S})|\\ &\leq
 n \tilde{C}_{0, 5}+n\epsilon |\lambda_{S, \epsilon }|,\\
\end{split}
\end{equation}
for all $0<\epsilon \leq 1$. Integrating the inequality (\ref{dd01}), we have
\begin{equation}\label{dd02}
\max_{x\in \tilde{M}}\epsilon |\log \frac{\det (\tilde{H}_{S, \epsilon }(1))}{ \det (H_{S, \epsilon }(0))}|\leq
 n \tilde{C}_{0, 5}+n\epsilon |\lambda_{S, \epsilon }|,
\end{equation}
for all $0<\epsilon \leq 1$. On the other hand, it is easy to check that there exists a uniform constant $\tilde{C}_{0, 6}$ such that:
\begin{equation}\label{100}
  \int_{\tilde{M} } |\sqrt{-1}\Lambda_{\omega_{\epsilon} }F_{H_{S, \epsilon }(0)}-\lambda_{S, \epsilon }\Id_{S}|_{H_{S, \epsilon }(0)}\frac{\omega_{\epsilon}^{n}}{n!}\leq \tilde{C}_{0, 6},
\end{equation}
for all $0<\epsilon \leq 1$. From the heat flow (\ref{SSS11}), one can check that
\begin{equation}
(\Delta_{\omega }-\frac{\partial }{\partial t})|\sqrt{-1}\Lambda_{\omega_{\epsilon} }F_{\tilde{H}_{S, \epsilon }(t)}-\lambda_{S, \epsilon }\Id_{S}|_{\tilde{H}_{S, \epsilon }(t)}\geq 0.
\end{equation}
By the maximum principle and the uniform upper bounds on the heat kernels (\ref{kernel01}), we have
\begin{equation}\label{curvature01}
\max_{x\in \tilde{M}}|\sqrt{-1}\Lambda_{\omega_{\epsilon} }F_{\tilde{H}_{S, \epsilon }(t)}-\lambda_{S, \epsilon }\Id_{S}|_{\tilde{H}_{S, \epsilon }(t)}\leq C_{K}\tilde{C}_{0, 6} (1+t^{-n}).
\end{equation}
This implies that there is a uniform constant $\tilde{C}_{0, 7}$ such that
\begin{equation}\label{dd03}
\max_{x\in \tilde{M}}\log \{\tr (\tilde{H}_{S, \epsilon }^{-1}(1)K_{S, \epsilon })+ \tr (K_{S, \epsilon }^{-1}\tilde{H}_{S, \epsilon }(1))\}\leq
 \tilde{C}_{0, 7}.
\end{equation}
From (\ref{dd02}) and (\ref{dd03}), we know that there exists a constant $\tilde{C}_{0, 8}$ such that
\begin{equation}\label{dd04}
\max_{x\in \tilde{M}}\epsilon \cdot |\log \frac{\det (K_{S, \epsilon })}{ \det (H_{S, \epsilon }(0))}|\leq \tilde{C}_{0, 8}.
\end{equation}

Let $\xi $ be the defining section of the line bundle $L$ with respect to the exceptional divisor $D\subset \pi^{-1}\Sigma_{alg}$. By the definition (\ref{Def05}) of the K\"ahler metric $\eta $, we see that
 \begin{equation}\label{Def06}
   \eta =\pi^{\ast }\omega -\delta \cdot \sqrt{-1}\partial \overline{\partial }\log |\xi |_{H_{L}}^{2},
   \end{equation}
on $\tilde{M}\setminus D$. Using (\ref{SQ}) and (\ref{80}), we get
\begin{equation}\label{701}
\begin{split}
&\int_{\tilde{M} } \tr (\varrho_{S, \epsilon }+\tilde{C}_{0, 2}\Id_{S}) \cdot \tr(\sqrt{-1}\Lambda_{\omega_{\epsilon} }F_{H_{S, \epsilon }(0) })\frac{\omega_{\epsilon}^{n}}{n!}\\
= & \int_{\tilde{M} } \tr (\varrho_{S, \epsilon }+\tilde{C}_{0, 2}\Id_{S}) \cdot \tr(\sqrt{-1}F_{H_{S, \epsilon }(0) })\wedge \frac{\omega_{\epsilon}^{n-1}}{(n-1)!}\\
\geq & -\tilde{C}_{0, 5}\int_{\tilde{M} } \tr (\varrho_{S, \epsilon }+\tilde{C}_{0, 2}\Id_{S}) \cdot \eta \wedge \frac{\omega_{\epsilon}^{n-1}}{(n-1)!}\\
= & -\tilde{C}_{0, 5}\int_{\tilde{M} } \tr (\varrho_{S, \epsilon }+\tilde{C}_{0, 2}\Id_{S}) \cdot (\pi^{\ast }\omega -\delta \cdot \sqrt{-1}\partial \overline{\partial }\log |\xi |_{H_{L}}^{2}) \wedge \frac{\omega_{\epsilon}^{n-1}}{(n-1)!}\\
\end{split}
\end{equation}

Near the divisor $D$, there is  a complex coordinate system $(U, (\tilde{z}^{1}, \cdots , \tilde{z}^{n}))$ such that $D\cap U=\{\tilde{z}^{n}\}$. So we can write locally:
\begin{equation}
|\xi |_{H_{L}}^{2}=\phi |\tilde{z}^{n}|^{2m},
\end{equation}
where $\phi $ is a nowhere vanishing smooth function. Setting $\beta_{1}=\sqrt{-1}\overline{\partial }\log |\xi |_{H_{L}}^{2}$, we see that $\beta_{1}$ is $1$-form with $L^{1}_{loc}$ coefficients on the entire $\tilde{M}$. By the residue formula, we have
\begin{equation}
\begin{split}
&\int_{\tilde{M} } \tr (\varrho_{S, \epsilon }) \cdot  \sqrt{-1}\partial \overline{\partial }\log |\xi |_{H_{L}}^{2} \wedge \frac{\omega_{\epsilon}^{n-1}}{(n-1)!}=\mathcal{T}_{d\beta_{1} }(\tr (\varrho_{S, \epsilon }) \frac{\omega_{\epsilon}^{n-1}}{(n-1)!})\\
= & d\mathcal{T}_{\beta_{1} }(\tr (\varrho_{S, \epsilon }) \frac{\omega_{\epsilon}^{n-1}}{(n-1)!}) +2m\pi \int_{D}\tr (\varrho_{S, \epsilon }) \frac{\omega_{\epsilon}^{n-1}}{(n-1)!}\\
= & \mathcal{T}_{\beta_{1} }(d\tr (\varrho_{S, \epsilon }) \wedge \frac{\omega_{\epsilon}^{n-1}}{(n-1)!}) +2m\pi \int_{D}\tr (\varrho_{S, \epsilon }) \frac{\omega_{\epsilon}^{n-1}}{(n-1)!},\\
\end{split}
\end{equation}
where $\mathcal{T}_{\beta }$ stands for the current with respect to the form $\beta$. Since $\tr(\varrho_{S, \epsilon })$ is a smooth function on $\tilde{M}$, $D\subset \pi^{-1}\Sigma_{alg}$ and $\Sigma_{alg} \subset M$ is a subset of complex codimension greater than or equal to $2$, we have
\begin{equation}
\int_{D}\tr (\varrho_{S, \epsilon }) \frac{(\pi^{\ast}\omega)^{n-1}}{(n-1)!} =0.
\end{equation}
By the estimate (\ref{dd04}), we obtain
\begin{equation}\label{702}
\begin{split}
& \int_{D}\tr (\varrho_{S, \epsilon }) \frac{\omega_{\epsilon}^{n-1}}{(n-1)!}\\
= &  \frac{1}{(n-1)!}\int_{D}\tr (\varrho_{S, \epsilon }) ((\pi^{\ast } \omega )^{n-1}+\epsilon (\Sigma_{i=1}^{n-1}C_{i}^{n-1}\epsilon^{i-1}\eta^{i}\wedge (\pi^{\ast }\omega )^{n-1-i})\\
= &  \frac{1}{(n-1)!}\int_{D}\epsilon \tr (\varrho_{S, \epsilon }) (\Sigma_{i=1}^{n-1}C_{i}^{n-1}\epsilon^{i-1}\eta^{i}\wedge (\pi^{\ast }\omega )^{n-1-i})\\
\geq & -\tilde{C}_{0, 9}\\
 \end{split}
\end{equation}
for all $0<\epsilon \leq 1$, where $\tilde{C}_{0, 9}$ is a positive constant.
Set $\beta_{2}=\sqrt{-1}\log |\xi |_{H_{L}}^{2}$, which is  $L^{1}_{loc}$  on the entire $\tilde{M}$. Using the residue formula again, we have
\begin{equation}\label{703}
\begin{split}
& \mathcal{T}_{\beta_{1} }(d\tr (\varrho_{S, \epsilon }) \wedge \frac{\omega_{\epsilon}^{n-1}}{(n-1)!})\\
= & \int_{\tilde{M} } \sqrt{-1}\overline{\partial }\log |\xi |_{H_{L}}^{2}\wedge \partial \tr (\varrho_{S, \epsilon }) \wedge \frac{\omega_{\epsilon}^{n-1}}{(n-1)!}\\
= &  \mathcal{T}_{d\beta_{2} }(\partial \tr (\varrho_{S, \epsilon }) \wedge \frac{\omega_{\epsilon}^{n-1}}{(n-1)!})\\
= &  - d\mathcal{T}_{\beta_{2} }(\partial \tr (\varrho_{S, \epsilon }) \wedge \frac{\omega_{\epsilon}^{n-1}}{(n-1)!}) -2m\pi \int_{D}\partial \tr (\varrho_{S, \epsilon }) \wedge \frac{\omega_{\epsilon}^{n-1}}{(n-1)!}\\
= &  \int_{\tilde{M} }\log |\xi |_{H_{L}}^{2}\sqrt{-1}\partial \overline{\partial } \tr (\varrho_{S, \epsilon }) \wedge \frac{\omega_{\epsilon}^{n-1}}{(n-1)!} -2m\pi \int_{D}\partial \tr (\varrho_{S, \epsilon }) \wedge \frac{\omega_{\epsilon}^{n-1}}{(n-1)!}.\\
\end{split}
\end{equation}
Since $D$ is a subset of complex codimension $1$, we see
\begin{equation}\label{704}
\int_{D}\partial \tr (\varrho_{S, \epsilon }) \wedge \frac{\omega_{\epsilon}^{n-1}}{(n-1)!}=0.
\end{equation}
On the other hand, it is easy to check that:
\begin{equation}\label{705}
\begin{split}
& \int_{\tilde{M} }\log |\xi |_{H_{L}}^{2}\sqrt{-1}\partial \overline{\partial } \tr (\varrho_{S, \epsilon }) \wedge \frac{\omega_{\epsilon}^{n-1}}{(n-1)!}\\
= & \int_{\tilde{M} }\log |\xi |_{H_{L}}^{2}\sqrt{-1}\partial \overline{\partial } \log \frac{\det (K_{S, \epsilon })}{ \det (H_{S, \epsilon }(0))} \wedge \frac{\omega_{\epsilon}^{n-1}}{(n-1)!}\\
= & \int_{\tilde{M} }\log |\xi |_{H_{L}}^{2}\sqrt{-1}\tr (F_{H_{S, \epsilon }(0)}-F_{K_{S, \epsilon }}) \wedge \frac{\omega_{\epsilon}^{n-1}}{(n-1)!}\\
= & \int_{\tilde{M} }\log |\xi |_{H_{L}}^{2}\{\sqrt{-1}\tr (F_{H_{S, \epsilon }(0)}) \wedge \frac{\omega_{\epsilon}^{n-1}}{(n-1)!}-\rank(S)\lambda_{S, \epsilon } \frac{\omega_{\epsilon}^{n}}{n!}\}\\
\geq &-\tilde{C}_{0, 10}\\
\end{split}
\end{equation}
for all $0<\epsilon \leq 1$, where $\tilde{C}_{0, 10}$ is a positive constant. Due to (\ref{701}), (\ref{702}), (\ref{703}), (\ref{704}) and (\ref{705}), there is a constant $\tilde{C}_{0, 11}$ such that
\begin{equation}\label{S70}
\int_{\tilde{M} } \tr (\varrho_{S, \epsilon }) \cdot \tr(\sqrt{-1}\Lambda_{\omega_{\epsilon} }F_{H_{S, \epsilon }(0) })\frac{\omega_{\epsilon}^{n}}{n!}\geq -\tilde{C}_{0, 11}
\end{equation}
for all $0<\epsilon \leq 1$. Then (\ref{700}) and (\ref{S70}) imply (\ref{MS}).

Applying the Gauss-Codazzi equation, we know
\begin{equation}\label{QQ70}
\tr(\sqrt{-1}\Lambda_{\omega_{\epsilon} }F_{H_{S, \epsilon }(0) }) +\tr (\sqrt{-1}\Lambda_{\omega_{\epsilon} }F_{H_{Q, \epsilon }(0) }) =\tr(\sqrt{-1}\Lambda_{\omega_{\epsilon} }F_{\pi^{\ast }H_{0} }),
\end{equation}
and
\begin{equation}\label{QQQQQ}
\begin{split}
&\Delta_{\epsilon} \{\tr (\varrho_{S, \epsilon }) + \tr (\varrho_{Q, \epsilon })\}\\
= & \tr(\sqrt{-1}\Lambda_{\omega_{\epsilon} }F_{H_{S, \epsilon }(0) }-\lambda_{S, \epsilon }\Id_{S})\\
  & +\tr(\sqrt{-1}\Lambda_{\omega_{\epsilon} }F_{H_{Q, \epsilon }(0) }-\lambda_{Q, \epsilon }\Id_{Q})\\
= & \tr(\sqrt{-1}\Lambda_{\omega_{\epsilon} }F_{\pi^{\ast }H_{0} }-\lambda_{\tilde{\mathcal{E}}, \epsilon }\Id_{\tilde{\mathcal{E}}}).
\end{split}
\end{equation}
By  (\ref{QQQQQ}), (\ref{bbin}), (\ref{det20}) and the uniform lower  bounds of the Green functions (\ref{green}), there is a constant $\tilde{C}_{0, 12}$ such that
\begin{equation}\label{QQQQQ1}
\max_{x\in \tilde{M}}|\tr (\varrho_{S, \epsilon }) + \tr (\varrho_{Q, \epsilon })|\leq \tilde{C}_{0, 12}.
\end{equation}
Then (\ref{q700}), (\ref{QQ70}), (\ref{QQQQQ1}) and (\ref{S70}) imply (\ref{MQ}).

\hfill $\Box$ \\

\medskip

\begin{proposition}  Along the heat flow (\ref{SSS1}), we have
\begin{equation}\label{C03}
\tr (\sqrt{-1}\Lambda_{\omega_{\epsilon} }(F_{H_{S, \epsilon }}-\gamma_{\epsilon} \wedge \gamma_{\epsilon }^{\ast})-\lambda_{S, \epsilon }\Id_{S})\leq \upsilon_{\epsilon} (t),
\end{equation}
for all $0<\epsilon \leq 1$, where $\upsilon_{\epsilon} (t)\geq 0$ satisfies: $\int_{1}^{+\infty} \upsilon (t) dt <\tilde{C}_{0, 13}$, and $\tilde{C}_{0, 13}$ is a uniform constant.
\end{proposition}

\medskip

{\bf Proof. }
 In \cite{D1}, for any exact sequence of holomorphic bundles $0\rightarrow S\rightarrow \tilde{\mathcal{E}} \rightarrow Q\rightarrow 0$, Donaldson has proved that
 \begin{equation}
 \begin{split}
  \mathcal{M} _{\tilde{\mathcal{E}}, \epsilon}^{0}(H_{0}, H(t))= &  \mathcal{M} _{S, \epsilon}^{0}(H_{S, 0}, H_{S})+ \mathcal{M} _{Q, \epsilon}^{0}(H_{Q, 0}, H_{Q})\\
 &+\|\gamma (t)\|_{L^{2}}^{2}-\|\gamma (0)\|_{L^{2}}^{2}.\\
 \end{split}
 \end{equation}
By the equations (\ref{S0}) and (\ref{Q0}), we have
\begin{equation}
\begin{split}
& \int_{\tilde{M}} \log \det (H_{S, \epsilon }^{-1}(0)H_{S, \epsilon }(t))\frac{\omega_{\epsilon}^{n}}{n!}\\=& \int_{0}^{t}\frac{\partial }{\partial l}\int_{\tilde{M}} \log \det (H_{S, \epsilon }^{-1}(0)H_{S, \epsilon }(l))\frac{\omega_{\epsilon}^{n}}{n!} dl\\
=& -2 \int_{0}^{t}\int_{\tilde{M}} \tr (\sqrt{-1}\Lambda_{\omega_{\epsilon } }(F_{H_{S, \epsilon }}-\gamma_{\epsilon} \wedge \gamma_{\epsilon }^{\ast})-\lambda_{\tilde{\mathcal{E}}}\Id_{S})\frac{\omega_{\epsilon}^{n}}{n!} dl\\
=& -2\int_{0}^{t}\int_{\tilde{M}}|\gamma_{\epsilon}(l) |_{H(l)}^{2} \frac{\omega_{\epsilon}^{n}}{n!} dl -2(\lambda_{S, \epsilon }-\lambda_{\tilde{\mathcal{E}}, \epsilon })\rank(S)t\Vol (\tilde{M}, \omega_{\epsilon } ),\\
\end{split}
\end{equation}
and
\begin{equation}
\int_{\tilde{M}} \log \det (H_{Q, \epsilon }^{-1}(0)H_{Q, \epsilon }(t))\frac{\omega_{\epsilon}^{n}}{n!}=2\int_{0}^{t}\int_{\tilde{M}}|\gamma_{\epsilon}(l) |_{H(l)}^{2} dl -2(\lambda_{Q_{\epsilon}}-\lambda_{\tilde{\mathcal{E}}})\rank(Q)t \Vol (\tilde{M}, \omega_{\epsilon } ) .\\
\end{equation}
Then
\begin{equation}
\begin{split}
& \mathcal{M} _{\tilde{\mathcal{E}}, \epsilon}(H_{0}, H) = \mathcal{M} _{S, \epsilon}(H_{S, \epsilon }(0) , H_{S, \epsilon }(t))+\mathcal{M} _{Q, \epsilon}(H_{Q, \epsilon }(0) , H_{Q, \epsilon }(t))\\
&+\|\gamma_{\epsilon } (t)\|_{L^{2}(\omega_{\epsilon })}^{2}-\|\gamma (0)\|_{L^{2}(\omega_{\epsilon })}^{2}\\
&-2(\lambda_{S, \epsilon }-\lambda_{\tilde{\mathcal{E}}, \epsilon })^{2}\rank(S)t\Vol (\tilde{M}, \omega_{\epsilon } ) -2(\lambda_{Q, \epsilon }-\lambda_{\tilde{\mathcal{E}}, \epsilon })^{2}\rank(Q)t\Vol (\tilde{M}, \omega_{\epsilon } ) \\
&-2 (\lambda_{S, \epsilon }-\lambda_{Q, \epsilon })\int_{0}^{t}\int_{\tilde{M}}|\gamma_{\epsilon}(l) |_{H(l)}^{2} \frac{\omega_{\epsilon}^{n}}{n!} dl.\\
\end{split}
\end{equation}
Furthermore, from (\ref{M0}), (\ref{C01}) and the Gauss-Codazzi equation (\ref{GC}), it follows that
\begin{equation}\label{sub1}
\begin{split}
&\mathcal{M} _{S, \epsilon}(H_{S, \epsilon }(0) , H_{S, \epsilon }(t))+\mathcal{M} _{Q, \epsilon}(H_{Q, \epsilon }(0) , H_{Q, \epsilon }(t))=-\|\gamma_{\epsilon } (t)\|_{L^{2}(\omega_{\epsilon })}^{2}+\|\gamma (0)\|_{L^{2}(\omega_{\epsilon })}^{2}\\
& -2 (\lambda_{S, \epsilon }-\lambda_{Q, \epsilon })\int_{0}^{t}\int_{\tilde{M}}|\gamma_{\epsilon}(l) |_{H(l)}^{2} \frac{\omega_{\epsilon}^{n}}{n!} dl -4\int_{0}^{t}\|\partial_{H_{S, \epsilon }}\gamma_{\epsilon } +\gamma_{\epsilon } \partial_{H_{Q, \epsilon }} \|_{L^{2}(\omega_{\epsilon})}^{2}(l)dl\\
& -2\int_{0}^{t}(\|\sqrt{-1}\Lambda_{\omega_{\epsilon} }(F_{H_{S, \epsilon }}-\gamma_{\epsilon} \wedge \gamma_{\epsilon}^{\ast })-\lambda_{S, \epsilon}\Id_{S}\|_{L^{2}(\omega_{\epsilon})}^{2})(l) dl
\\ & -2\int_{0}^{t}(\|\sqrt{-1}\Lambda_{\omega_{\epsilon} }(F_{H_{Q, \epsilon }}-\gamma_{\epsilon}^{\ast} \wedge \gamma_{\epsilon})-\lambda_{Q, \epsilon}\Id_{Q}\|_{L^{2}(\omega_{\epsilon})}^{2})(l) dl.\\
\end{split}
\end{equation}

Now, we set:
\begin{equation}
\begin{split}
\tilde{\upsilon }(\epsilon , t)=& 2\|\partial_{H_{S, \epsilon }}\gamma_{\epsilon } +\gamma_{\epsilon } \partial_{H_{Q, \epsilon }} \|_{L^{2}(\omega_{\epsilon})}^{2}(t)+ (\lambda_{S, \epsilon }-\lambda_{Q, \epsilon })\|\gamma_{\epsilon } \|_{L^{2}(\omega_{\epsilon})}^{2}(t)\\
&+\|\sqrt{-1}\Lambda_{\omega_{\epsilon} }(F_{H_{S, \epsilon }}-\gamma_{\epsilon} \wedge \gamma_{\epsilon}^{\ast })-\lambda_{S, \epsilon}\Id_{S}\|_{L^{2}(\omega_{\epsilon})}^{2}(t)\\
 &+ \|\sqrt{-1}\Lambda_{\omega_{\epsilon} }(F_{H_{Q, \epsilon }}-\gamma_{\epsilon}^{\ast} \wedge \gamma_{\epsilon})-\lambda_{Q, \epsilon}\Id_{Q}\|_{L^{2}(\omega_{\epsilon})}^{2}(t).\\
\end{split}
\end{equation}
Combining the uniform lower bounds of the Donaldson's functional (\ref{MS}), (\ref{MQ}) and the above equality (\ref{sub1}), we obtain
\begin{equation}\label{sub2}
\int_{0}^{+\infty} \tilde{\upsilon }(\epsilon , t) dt \leq  \tilde{C}_{0, 14} <+\infty ,
\end{equation}
for all $0<\epsilon \leq 1$, where $\tilde{C}_{0, 14}$ is a uniform positive constant.

On the other hand, along the heat flow (\ref{SSS1}) on $\tilde{\mathcal{E}}$, we have
\begin{equation}
(\Delta_{\epsilon } -\frac{\partial }{\partial t})|\sqrt{-1}\Lambda_{\omega_{\epsilon} }F_{H_{\epsilon}(t)}-\lambda_{\tilde{\mathcal{E}}, \epsilon } \Id_{\tilde{\mathcal{E}}}|_{H_{\epsilon}(t)}^{2}\geq 0.
\end{equation}
By the uniform  Sobolev inequality (\ref{Sobolev}) for $(\tilde{M}, \omega_{\epsilon})$ and the estimate of the heat kernel $\chi (x, y, t)$ by Cheng and Li ($(2.9)$ in \cite{CL}, or Theorem 3.2 in \cite{BM}), there is a positive constant $\tilde{C}_{0, 15}$ such that
\begin{equation}
0<\mathcal{K}_{\epsilon}(x , y, t)\leq \frac{1}{\Vol (\tilde{M}, \omega_{\epsilon } )} +\tilde{C}_{0, 15}t^{-n},
\end{equation}
for $0<\epsilon \leq 1$.
By the maximum principle,  we have
\begin{equation}\label{heatkernel}
\begin{split}
&|\sqrt{-1}\Lambda_{\omega_{\epsilon} }F_{H_{\epsilon}(t)}-\lambda_{\tilde{\mathcal{E}}, \epsilon } \Id_{\tilde{\mathcal{E}}}|_{H_{\epsilon}(t)}^{2}(x, t+s)\\
\leq & \int_{\tilde{M}}\mathcal{K}_{\epsilon}(x , y, s)|\sqrt{-1}\Lambda_{\omega_{\epsilon} }F_{H_{\epsilon}(t)}-\lambda_{\tilde{\mathcal{E}}, \epsilon } \Id_{\tilde{\mathcal{E}}}|_{H_{\epsilon}(t)}^{2}(y, t) \frac{\omega_{\epsilon}^{n}}{n!}(y)\\
\leq & (\frac{1}{\Vol (\tilde{M}, \omega_{\epsilon } )} +\tilde{C}_{0, 15}s^{-n})\int_{\tilde{M}}|\sqrt{-1}\Lambda_{\omega_{\epsilon} }F_{H_{\epsilon}(t)}-\lambda_{\tilde{\mathcal{E}}, \epsilon } \Id_{\tilde{\mathcal{E}}}|_{H_{\epsilon}(t)}^{2}(y, t) \frac{\omega_{\epsilon}^{n}}{n!}(y),\\
\end{split}
\end{equation}
for any $t>0$, $s>0$ and $0<\epsilon \leq 1$.

Using the Gauss-Codazzi equation (\ref{GC}) again, and combining the formulas (\ref{C01}) and (\ref{heatkernel}), we have
\begin{equation}\label{xi}
\begin{split}
&2(\lambda_{S, \epsilon }-\lambda_{Q, \epsilon })\tr (\sqrt{-1}\Lambda_{\omega_{\epsilon} }(F_{H_{S, \epsilon }}-\gamma_{\epsilon} \wedge \gamma_{\epsilon}^{\ast })-\lambda_{S, \epsilon}\Id_{S})(t)\\
\leq &|\sqrt{-1}\Lambda_{\omega_{\epsilon} }F_{H_{\epsilon}(t)}-\lambda_{\tilde{\mathcal{E}}, \epsilon } \Id_{\tilde{\mathcal{E}}}|_{H_{\epsilon}(t)}^{2}\\ &-(\lambda_{S, \epsilon }-\lambda_{\tilde{\mathcal{E}}, \epsilon })^{2}\rank(S)-(\lambda_{Q, \epsilon }-\lambda_{\tilde{\mathcal{E}}, \epsilon })^{2}\rank(Q)\\
& - 2(\lambda_{Q, \epsilon }-\lambda_{\mathcal{E}, \epsilon }) \tr (\sqrt{-1}\Lambda_{\omega_{\epsilon } }F_{H_{\epsilon} (t)} -\lambda_{\mathcal{E}, \epsilon } \Id_{\mathcal{E}})\\
\leq & (\frac{1}{\Vol (\tilde{M}, \omega_{\epsilon } )} +\tilde{C}_{0, 15}(\frac{t}{2})^{-n})\int_{\tilde{M}}|\sqrt{-1}\Lambda_{\omega_{\epsilon} }F_{H_{\epsilon}(\frac{t}{2})}-\lambda_{\tilde{\mathcal{E}}, \epsilon } \Id_{\tilde{\mathcal{E}}}|_{H_{\epsilon}(\frac{t}{2})}^{2} \\ &-(\lambda_{S, \epsilon }-\lambda_{\tilde{\mathcal{E}}, \epsilon })^{2}\rank(S)-(\lambda_{Q, \epsilon }-\lambda_{\tilde{\mathcal{E}}, \epsilon })^{2}\rank(Q)+ 2(\lambda_{\mathcal{E}, \epsilon }-\lambda_{Q, \epsilon })|f_{\epsilon }(t)|\\
\leq & (\frac{1}{\Vol (\tilde{M}, \omega_{\epsilon } )}+\tilde{C}_{0, 15}(\frac{t}{2})^{-n})\tilde{\upsilon }(\epsilon, \frac{t}{2})+2(\lambda_{\mathcal{E}, \epsilon }-\lambda_{Q, \epsilon })|f_{\epsilon }(t)|\\
& +\Vol (\tilde{M}, \omega_{\epsilon } )\tilde{C}_{0, 15}(\frac{t}{2})^{-n}\{(\lambda_{S, \epsilon }-\lambda_{\mathcal{E}, \epsilon })^{2}\rank(S)+(\lambda_{Q, \epsilon } -\lambda_{\mathcal{E}, \epsilon})^{2}\rank(Q)\}. \\
\end{split}
\end{equation}
Set
\begin{equation}
\begin{split}
\upsilon_{\epsilon} (t)=&(2(\lambda_{S, \epsilon }-\lambda_{Q, \epsilon }))^{-1}\{ (\frac{1}{\Vol (\tilde{M}, \omega_{\epsilon } )}+\tilde{C}_{0, 15}(\frac{t}{2})^{-n})\tilde{\upsilon}(\epsilon, \frac{t}{2})+2(\lambda_{\mathcal{E}, \epsilon }-\lambda_{Q, \epsilon })\max_{x\in \tilde{M}}|f_{\epsilon }(x, t)|\\
& +\Vol (\tilde{M}, \omega_{\epsilon } )\tilde{C}_{0, 15}(\frac{t}{2})^{-n}\{(\lambda_{S, \epsilon }-\lambda_{\mathcal{E}, \epsilon })^{2}\rank(S)+(\lambda_{Q, \epsilon} -\lambda_{\mathcal{E}, \epsilon })^{2}\rank(Q)\}\}.\\
\end{split}
\end{equation}
By the formulas (\ref{sub2}), (\ref{xi}), (\ref{FF1}) and note that $n>1$, we see that $\upsilon (t)$
is the function which we need.

\hfill $\Box$ \\

\medskip

Using the above proposition, we can obtain a uniform local $C^{0}$-bound on the rescaled metrics $\tilde{H}_{S}(t)=e^{2(\lambda_{S}-\lambda_{E})t}H_{S}(t)$ and $\tilde{H}_{Q}(t)=e^{2(\lambda_{Q}-\lambda_{E})t}H_{Q}(t)$.

\medskip

\begin{theorem}  Let $H(t)$ be the solution of the Hermitian-Yang-Mills flow (\ref{Dflow}) on the holomorphic bundle $\mathcal{E}$ with initial metric $H_{0}$,  $H_{\mathcal{S}}(t)$ and $H_{\mathcal{Q}}(t)$ be the induced Hermitian metrics on $\mathcal{S}|_{M\setminus \Sigma_{alg}}$ and $\mathcal{Q}|_{M\setminus \Sigma_{alg}}$.  Set $\hat{h}_{\mathcal{S}}(t)=e^{2(\lambda_{\mathcal{S}}-\lambda_{\mathcal{E}})t}H_{\mathcal{S}}(0)^{-1}H_{\mathcal{S}}(t)$ and $\hat{h}_{\mathcal{Q}}(t)=e^{2(\lambda_{\mathcal{Q}}-\lambda_{\mathcal{E}})t}H_{\mathcal{Q}}(0)^{-1}H_{\mathcal{Q}}(t)$. Then for any compact subset $U\subset M\setminus \Sigma_{alg}$ there exists a uniform constant $\tilde{C}_{U}$ such that
\begin{equation}\label{C0}
\tr \hat{h}_{\mathcal{S}}(x, t) +\tr \hat{h}_{\mathcal{S}}^{-1} (x, t) + \tr \hat{h}_{\mathcal{Q}}(x, t)+ \tr \hat{h}_{\mathcal{Q}}^{-1}(x, t)\leq \tilde{C}_{U},
\end{equation}
for all $t\geq 0$ and $x\in U$.
\end{theorem}

\medskip

{\bf Proof. } In the following, we denote:
\begin{equation}
B_{\Sigma }(\delta )=\{x\in M|d(x, \Sigma_{alg } )<\delta \},
\end{equation}
where $d$ is the distance function with respect to the K\"ahler metric $\omega $. We can choose  a small number $\delta_{0}$ such that $U\subset M\setminus B_{\Sigma }(\delta_{0} )$. Since $\pi^{\ast} \omega $ is positive on $\tilde{M}\setminus \pi^{-1}\Sigma_{alg}$, there is a constant $C_{\delta_{0}^{-1}}$ such that
\begin{equation}
C_{\delta_{0}^{-1}}^{-1}\eta \leq  \pi^{\ast} \omega \leq C_{\delta_{0}^{-1}}\eta
\end{equation}
on $\pi^{-1}(M\setminus B_{\Sigma }(\frac{1}{2}\delta_{0} ))$. Noting that the metrics $H_{S, \epsilon }(0)$ and $H_{Q, \epsilon }(0)$ are independent of $\epsilon$, one checks that
\begin{equation}\label{FS1}
 -C_{U, 1} \Id_{S} \leq \sqrt{-1}\Lambda_{\omega_{\epsilon}}  F_{H_{S, \epsilon}(0)} \leq C_{U, 1} \Id_{S},
\end{equation}
and
\begin{equation}\label{FQ1}
-C_{U, 1} \Id_{Q} \leq  \sqrt{-1}\Lambda_{\omega_{\epsilon}}  F_{H_{Q, \epsilon}(0)} \leq C_{U, 1} \Id_{Q}
\end{equation}
on $\pi^{-1}(M\setminus B_{\Sigma }(\frac{1}{2}\delta_{0} ))$, where $C_{U, 1}$ is a uniform constant independent of $\epsilon$. Then there is a uniform constant $C_{U, 2}$ such that
\begin{equation}\label{TTTT10011}
\Delta_{\epsilon} \log (\tr \tilde{h}_{S, \epsilon}(0)+ \tr \tilde{h}_{S, \epsilon}^{-1}(0)) \geq   -C_{U, 2}
\end{equation}
and
\begin{equation}\label{TTTT10021}
\Delta_{\epsilon} \log (\tr \tilde{h}_{Q, \epsilon}(0)+ \tr \tilde{h}_{Q, \epsilon}^{-1}(0)) \geq   -C_{U, 2}
\end{equation}
on $\pi^{-1}(M\setminus B_{\Sigma }(\frac{1}{2}\delta_{0} ))$, for all $0<\epsilon \leq 1$. By  the uniform Sobolev inequality (\ref{Sobolev}), (\ref{TTTT10011}), (\ref{TTTT10021}) and Moser's iteration,  we have the following mean-value inequalities, i.e.  there exists a uniform constant $C_{U, 3}$ such that
\begin{equation}
\sup_{\pi^{-1}(M\setminus B_{\Sigma }(\delta_{0} ))} \log (\tr \tilde{h}_{S, \epsilon}(0)+ \tr \tilde{h}_{S, \epsilon}^{-1}(0)) \leq C_{U, 3} \int_{\tilde{M}}\log \{\tr(\tilde{h}_{S, \epsilon}(0))+\tr(\tilde{h}_{S, \epsilon}^{-1}(0))\}\frac{\omega_{\epsilon}^{n}}{n!}
\end{equation}
and
\begin{equation}
\sup_{\pi^{-1}(M\setminus B_{\Sigma }(\delta_{0} ))} \log (\tr \tilde{h}_{Q, \epsilon}(0)+ \tr \tilde{h}_{Q, \epsilon}^{-1}(0)) \leq C_{U, 3} \int_{\tilde{M}}\log \{\tr(\tilde{h}_{Q, \epsilon}(0))+\tr(\tilde{h}_{Q, \epsilon}^{-1}(0))\}\frac{\omega_{\epsilon}^{n}}{n!}.
\end{equation}
From (\ref{det20}), we know that there is a uniform constant $C_{U, 4}$ such that
\begin{equation}\label{CC00}
\sup_{\pi^{-1}(M\setminus B_{\Sigma }(\delta_{0} ))} \{\log (\tr \tilde{h}_{S, \epsilon}(0)+ \tr \tilde{h}_{S, \epsilon}^{-1}(0))+  \log (\tr \tilde{h}_{Q, \epsilon}(0)+ \tr \tilde{h}_{Q, \epsilon}^{-1}(0))\}\leq C_{U, 4},
\end{equation}
for all $0<\epsilon \leq 1$. Set
\[
\hat{h}_{S, \epsilon}(t)=e^{2(\lambda_{S, \epsilon }-\lambda_{\tilde{\mathcal{E}}, \epsilon })t}H_{S, \epsilon}(0)^{-1}H_{S, \epsilon}(t)
\]
and
\[
\hat{h}_{Q, \epsilon}(t)=e^{2(\lambda_{Q, \epsilon }-\lambda_{\tilde{\mathcal{E}}, \epsilon })t}H_{Q, \epsilon}(0)^{-1}H_{Q, \epsilon}(t).
\]
(\ref{C00}) and (\ref{CC00}) imply that
\begin{equation}\label{CCC00}
\sup_{(x, t)\in \pi^{-1}(M\setminus B_{\Sigma }(\delta_{0} ))\times [0, +\infty )}(\tr \hat{h}_{S, \epsilon }(x, t) + \tr \hat{h}_{Q, \epsilon }^{-1}(x, t))\leq C_{U, 5},
\end{equation}
for all $0<\epsilon \leq 1$, where $C_{U, 5}$ is a uniform constant independent of $\epsilon$.
Due to (\ref{S0}), it holds that
\begin{equation}
\frac{\partial }{\partial t}\log (\det \hat{h}_{S, \epsilon }^{-1})=2 \tr (\sqrt{-1}\Lambda_{\omega_{\epsilon} }(F_{H_{S, \epsilon }}-\gamma_{\epsilon } \wedge \gamma_{\epsilon }^{\ast})-\lambda_{S, \epsilon }\Id_{S}).
\end{equation}
By (\ref{bbin}), the Gauss-Codazzi equation (\ref{GC}) and the maximum principle, we have
\begin{equation}
\sup_{(x, t)\in \tilde{M}\times [0, +\infty )}|\sqrt{-1}\Lambda_{\omega_{\epsilon}}F_{H_{\epsilon}(t)}-\lambda_{\tilde{\mathcal{E}}, \epsilon }Id_{\tilde{\mathcal{E}}}|^{2}_{H_{\epsilon}(t)}\leq \overline{C}_{F}
\end{equation}
and
\begin{equation}
\sup_{(x, t)\in \tilde{M}\times [0, +\infty )}|\tr (\sqrt{-1}\Lambda_{\omega_{\epsilon} }(F_{H_{S, \epsilon }}-\gamma_{\epsilon } \wedge \gamma_{\epsilon }^{\ast})-\lambda_{S, \epsilon }\Id_{S})|\leq \overline{C}_{F},
\end{equation}
where $\overline{C}_{F}$ is a uniform constant independent of $\epsilon$.
Then (\ref{C03}) implies that there exists a uniform constant $C_{U, 6}$ such that
\begin{equation}\label{C05}
\sup_{(x, t)\in \tilde{M}\times [0, +\infty )}\log (\det \hat{h}_{S, \epsilon }^{-1}(x, t))\leq C_{U, 6},
\end{equation}
for all $0<\epsilon \leq 1$.
Combining (\ref{C04}), (\ref{FF1}) and (\ref{C05}), we see that
\begin{equation}
\det(\hat{h}_{S, \epsilon }( t))\det(\hat{h}_{Q, \epsilon }(t))=\det(H_{0}^{-1}H_{\epsilon}(t)),
\end{equation}
and then
\begin{equation}\label{C06}
\begin{split}
\log \det (\hat{h}_{Q, \epsilon }(x, t))= & \log \det (\hat{h}_{S, \epsilon }^{-1}(x, t))+\log \det(H_{0}^{-1}H_{\epsilon}(x, t))\\
\leq & C_{U, 7}\\
\end{split}
\end{equation}
for all $(x, t)\in \tilde{M}\times [0, +\infty )$ and $0<\epsilon \leq 1$, where $C_{U, 7}$ is a uniform constant independent of $\epsilon $.
 From (\ref{CCC00}), (\ref{C05}) and (\ref{C06}), it is easy to see that there exists a constant $C_{U, 8}$ such that
\begin{equation}\label{CCC00}
\sup_{(x, t)\in \pi^{-1}(M\setminus B_{\Sigma }(\delta_{0} ))\times [0, +\infty )}(\tr \hat{h}_{S, \epsilon }+\tr \hat{h}_{S, \epsilon }^{-1} +\tr \hat{h}_{Q, \epsilon }+ \tr \hat{h}_{Q, \epsilon }^{-1})(x, t)\leq C_{U, 8}
\end{equation}
for all $0<\epsilon \leq 1 $. Taking the limit $\epsilon \rightarrow 0$, we obtain (\ref{C0}).

\hfill $\Box$ \\

\medskip

\section{Uniform local $C^{1}$-estimate }
\setcounter{equation}{0}

Using the above local  $C^{0}$-estimate of the rescaled metrics $\tilde{H}_{S}(t)$ and $\tilde{H}_{Q}(t)$, we can control the  $L_{loc}^{\infty}$-norm of $|G(t)|_{H(t)}$. In fact, we have the following proposition.

\medskip

\begin{proposition}  Let $H(t)$ be the solution of the heat flow (\ref{Dflow}) with initial metric $H_{0}$ on $(E, \overline{\partial} _{E})$, and $G(t)\in \Gamma (\mathcal{S}\otimes \mathcal{Q}^{\ast}|_{M\setminus \Sigma_{alg}})$ be the section defined in (\ref{is2}). Then for any compact subset $U\subset M\setminus \Sigma_{alg}$, there exists a constant $C_{G, U}$ such that
\begin{equation}\label{CG}
\int_{U}|G (t)|_{H(t)}^{2} \frac{\omega^{n}}{n!}\leq C_{G, U}
\end{equation}
for any $t\geq 0$. Furthermore,  there exists a constant $\tilde{C}_{G, U}$ such that
\begin{equation}\label{CG2}
\sup_{(x, t)\in U\times [0, +\infty ) }|G (t)|_{H(t)}^{2}(x) \leq \tilde{C}_{G, U}.
\end{equation}
\end{proposition}

\medskip

{\bf Proof. } Let $\delta $ be small enough so that $U\subset M\setminus B_{\Sigma }(\delta )$, and  $\varphi$ be a  nonnegative cut-off function satisfying:
\begin{equation}
 \varphi (x)=\left \{\begin{split} &1,   \quad  x\in M\setminus B_{\Sigma }(\delta ),\\
&0,  \quad x\in  B_{\Sigma }(\frac{1}{2}\delta ).\\
\end{split}
\right.
\end{equation}
Direct calculations yield
\begin{equation}\label{G14}
\begin{split}
&\frac{\partial }{\partial t}\int_{M}\varphi ^{2}e^{-2(\lambda_{S}-\lambda_{Q})t} |G (t)|_{H(0)}^{2} \frac{\omega^{n}}{n!}\\
= & \int_{M} 2\varphi ^{2}e^{-2(\lambda_{S}-\lambda_{Q})t}Re\big\langle \frac{\partial G }{\partial t}, G(t)\big\rangle_{H(0)}\frac{\omega^{n}}{n!}\\
& -2(\lambda_{S}-\lambda_{Q})\int_{M}\varphi ^{2}e^{-2(\lambda_{S}-\lambda_{Q})t} |G (t)|_{H(0)}^{2} \frac{\omega^{n}}{n!}\\
\leq & 2\{\int_{M} \varphi ^{2}e^{-2(\lambda_{S}-\lambda_{Q})t}\big| \frac{\partial G }{\partial t}\big|_{H(0)}^{2}\frac{\omega^{n}}{n!}\}^{\frac{1}{2}}\{\int_{M}\varphi ^{2}e^{-2(\lambda_{S}-\lambda_{Q})t} |G (t)|_{H(0)}^{2} \frac{\omega^{n}}{n!}\}^{\frac{1}{2}}\\
& -2(\lambda_{S}-\lambda_{Q})\int_{M}\varphi ^{2}e^{-2(\lambda_{S}-\lambda_{Q})t} |G (t)|_{H(0)}^{2} \frac{\omega^{n}}{n!}\\
\leq & \frac{1}{2(\lambda_{S}-\lambda_{Q})}\int_{M} \varphi ^{2}e^{-2(\lambda_{S}-\lambda_{Q})t}\big| \frac{\partial G }{\partial t}\big|_{H(0)}^{2}\frac{\omega^{n}}{n!}.\\
\end{split}
\end{equation}
For any section $\theta \in \Gamma (\mathcal{S}\otimes \mathcal{Q}^{\ast})$, by the local $C^{0}$-estimate  (\ref{C0}), we know that there is a uniform constant $C_{\delta^{-1}, 1}$ such that
\begin{equation}\label{G15}
\begin{split}
& |\theta |_{H(0)}^{2} (x) = \tr (\theta\cdot H_{Q}^{-1}(0)\overline{\theta}^{T} H_{S}(0))(x)\\
\leq & C_{\delta^{-1}, 1}\tr (\theta \cdot e^{-2(\lambda_{Q}-\lambda_{E})t}H_{Q}^{-1}(t)\overline{\theta}^{T} e^{2(\lambda_{S}-\lambda_{E})t}H_{S}(t))(x)\\
= & e^{2(\lambda_{S}-\lambda_{Q})t}C_{\delta^{-1}, 1} \tr (\theta \cdot H_{Q}^{-1}(t)\overline{\theta}^{T} H_{S}(t))(x)\\
= & C_{\delta^{-1}, 1} e^{2(\lambda_{S}-\lambda_{Q})t} |\theta|_{H(t)}^{2}(x)\\
\end{split}
\end{equation}
and
\begin{equation}\label{G151}
|\theta |_{H(0)}^{2} (x)\geq C_{\delta^{-1}, 1}^{-1}e^{2(\lambda_{S}-\lambda_{Q})t} |\theta|_{H(t)}^{2}(x)\\
\end{equation}
for all $(x, t)\in M\setminus B_{\Sigma }(\frac{1}{2}\delta )\times [0, +\infty )$. Using (\ref{G14}), (\ref{G15}), (\ref{G151}) (\ref{G0}) and (\ref{sub2}), we get
\begin{equation}\label{G141}
\begin{split}
& \int_{M}\varphi ^{2} |G (t)|_{H(t)}^{2} \frac{\omega^{n}}{n!}\\
\leq & C_{\delta^{-1}, 1}\int_{M}\varphi ^{2}e^{-2(\lambda_{S}-\lambda_{Q})t} |G (t)|_{H(0)}^{2} \frac{\omega^{n}}{n!}\\
= & C_{\delta^{-1}, 1}\int_{0}^{t}\frac{\partial }{\partial l}\int_{M}\varphi ^{2}e^{-2(\lambda_{S}-\lambda_{Q})l} |G (l)|_{H(0)}^{2} \frac{\omega^{n}}{n!} dl\\
\leq &  C_{\delta^{-1}, 1}^{2}\frac{1}{2(\lambda_{S}-\lambda_{Q})}\int_{M} \varphi ^{2}| \Lambda_{\omega}(\partial_{H_{S}(t)}\gamma(t)+ \gamma(t)\partial_{H_{Q}(t)} )|_{H(t)}^{2}\frac{\omega^{n}}{n!}\\
\leq & C_{\delta^{-1}, 2},\\
\end{split}
\end{equation}
for all $t\in [0, +\infty )$. It is easy to see that (\ref{G141}) implies (\ref{CG}).

\medskip

Since $\overline{\partial }_{S\otimes Q^{\ast}}\gamma_{0} =0$, for any point $P\in U$, we have a domain $B_{P}(4R)$ and a local section  $G_{0}\in \Gamma (B_{P}(4R); S\otimes Q^{\ast})$ such that $\gamma_{0}=\overline{\partial }_{S\otimes Q^{\ast}} G_{0}$, where $R< \frac{1}{8}\delta $. By (\ref{G6}),
we have
\begin{equation}
(\Delta -\frac{\partial }{\partial t})|G+G_{0}|_{H(t)}^{2}\geq 0
\end{equation}
on $B_{P}(4R)$.
 Applying the parabolic mean value inequality (Theorem 14.5. in \cite{Pli2}, or \cite{LT}), we obtain
\begin{equation}
\sup_{B_{P}(\frac{R}{2})\times [t_{0}+\frac{R^{2}}{8},t_{0}+\frac{R^{2}}{4}] }|G+G_{0}|_{H(t)}^{2}\leq C_{1}\int_{\frac{R^{2}}{16}}^{\frac{R^{2}}{2}}\int_{B_{P}(R)}|G+G_{0}|_{H(s+t_{0})}^{2}dv_{\omega}ds
\end{equation}
for any $t_{0} >0$, where constant $C_{1}$ depends only on dim$M$, lower bound of Ricci curvature and $R^{-1}$. As in (\ref{G151}), we know $|G_{0}|_{H(t)}$ is uniformly bounded on $B_{P}(2R)\times [0, +\infty )$. From (\ref{CG}), we see
 $|G|_{H(t)}$ is also uniformly bounded on $B_{P}(\frac{R}{2})\times [0, +\infty )$. Since point $P$ is arbitrary and $U$ is compact,  we can obtain a uniform bound for $|G|_{H(t)}$ on $U\times [0, +\infty )$, i.e. the inequality (\ref{CG2}) is valid.

\hfill $\Box$ \\

Set
\begin{equation}\label{DTS1}
T_{S}(t)=D_{(H_{S}(t), \overline{\partial}_{S})}-D_{(H_{0, S}, \overline{\partial}_{S})}=h_{S}^{-1}\partial_{H_{0, S}}h_{S}=(\partial_{H_{S}(t)}h_{S})h_{S}^{-1}
\end{equation}
and
\begin{equation}\label{DTQ1}
T_{Q}(t)=D_{(H_{Q}(t), \overline{\partial}_{Q})}-D_{(H_{0, Q}, \overline{\partial}_{Q})}=h_{S}^{-1}\partial_{H_{0, Q}}h_{Q}=(\partial_{H_{Q}(t)}h_{Q})h_{Q}^{-1},
\end{equation}
where $h_{S}(t)=H_{0, S}^{-1}H_{S}(t)$ and $h_{Q}(t)=H_{0, Q}^{-1}H_{Q}(t)$. By the definition, we have
\begin{equation}
\begin{split}
\frac{\partial }{\partial t}T_{S}=& \frac{\partial }{\partial t}(H_{S}^{-1}(t)\partial H_{S}(t))\\
=& \partial_{H_{S}(t)}(H_{S}^{-1}(t)\frac{\partial H_{S}(t)}{\partial t}),\\
\end{split}
\end{equation}
and
\begin{equation}\label{TS0}
\begin{split}
\frac{\partial }{\partial t}|T_{S}|_{H_{S}}^{2}=& \frac{\partial }{\partial t}g^{k\bar{l}}\tr(T_{S}(\partial_{k})H_{S}^{-1}\overline{T_{S}(\partial_{l})^{T}}H_{S})\\
=& 2 Re \{ g^{k\bar{l}}\tr(\frac{\partial T_{S} }{\partial t}(\partial_{k})H_{S}^{-1}\overline{T_{S}(\partial_{l})^{T}}H_{S})\}\\
&-g^{k\bar{l}}\tr(T_{S}(\partial_{k})H_{S}^{-1}\frac{\partial H_{S}}{\partial t}H_{S}^{-1}\overline{T_{S}(\partial_{l})^{T}}H_{S})\\
&+g^{k\bar{l}}\tr(T_{S}(\partial_{k})H_{S}^{-1}\overline{T_{S}(\partial_{l})^{T}}H_{S}H_{S}^{-1}\frac{\partial H_{S}}{\partial t})\\
=& 2Re \langle\partial_{H_{S}}(H_{S}^{-1}\frac{\partial H_{S}}{\partial t}), T_{S}\rangle_{H_{S}(t)}\\
&+Re \langle[H_{S}^{-1}\frac{\partial H_{S}}{\partial t}, T_{S} ] , T_{S}\rangle_{H_{S}(t)}.\\
\end{split}
\end{equation}
Now we calculate the Laplacian of $|T_{S}|_{H_{S}}^{2}$. In the following we choose the normal complex coordinates $\{z^{1}, \cdots , z^{m}\}$ centered at the considering point, and denote $\frac{\partial}{\partial z^{i}}$ ($\frac{\partial}{\partial \overline{z}^{j}}$) by $\partial _{i}$ ($\overline{\partial }_{j}$) for simplicity. Direct  calculations yield
\begin{equation}
\begin{split}
\Delta |T_{S}|_{H_{S}}^{2} =& 2g^{i\bar{j}}\partial_{i}\overline{\partial}_{j} |T_{S}|_{H_{S}}^{2}\\
=& 2Re \{g^{i\bar{j}} \langle\nabla _{\partial_{i}}^{H_{S}(t)}\nabla _{\overline{\partial}_{j}}^{H_{S}(t)}T_{S} +\nabla _{\overline{\partial}_{j}}^{H(t)}\nabla _{\partial_{i}}^{H(t)}T_{S}, T_{S} \rangle_{H_{S}(t)}\}\\
&+2|\nabla ^{H_{S}(t)}T_{S} |_{H_{S}(t)}^{2},\\
\end{split}
\end{equation}
\begin{equation}
\begin{split}
\nabla _{\partial_{i}}^{H_{S}(t)}\nabla _{\overline{\partial}_{j}}^{H_{S}(t)}T_{S} (\partial_{k}) =& D _{\partial_{i}}^{H_{S}(t)}(\nabla _{\overline{\partial}_{j}}^{H_{S}(t)}T_{S}(\partial_{k}))-(\nabla _{\overline{\partial}_{j}}^{H_{S}(t)}T_{S})(\nabla_{\partial_{i}}\partial_{k})\\
=& D _{\partial_{i}}^{H_{S}(t)}(F_{H_{0, S}}(\partial_{k} , \overline{\partial}_{j})-F_{H_{S}}(\partial_{k} , \overline{\partial}_{j}))\\
=& D _{\partial_{i}}^{H_{S}(t)}(F_{H_{0, S}}(\partial_{k} , \overline{\partial}_{j}))- D _{\partial_{k}}^{H_{S}(t)}(F_{H_{S}}(\partial_{i} , \overline{\partial}_{j})),\\
\end{split}
\end{equation}
where we have used the equality $F_{H_{S}}=F_{H_{0, S}}+\overline{\partial }_{S}(T_{S})$, the Bianchi identity $(\nabla _{\partial_{i}}^{H_{S}(t)}F_{H_{S}})(\partial_{k} , \overline{\partial}_{j})=(\nabla _{\partial_{k}}^{H_{S}(t)}F_{H_{S}})(\partial_{i} , \overline{\partial}_{j})$, and $\nabla_{\partial_{i}}\partial_{k} =0$ at the considering point.
Furthermore,
\begin{equation}
\begin{split}
\nabla _{\overline{\partial}_{j}}^{H_{S}(t)}\nabla _{\partial_{i}}^{H_{S}(t)}T_{S} (\partial_{k}) =& \nabla _{\overline{\partial}_{j}}^{H_{S}(t)}(D _{\partial_{i}}^{H_{S}(t)}(T_{S}(\partial_{k}))-T_{S}(\nabla_{\partial_{i}}\partial_{k}))\\
=& \nabla _{\overline{\partial}_{j}}^{H_{S}(t)}D _{\partial_{i}}^{H_{S}(t)}(T_{S}(\partial_{k}))-T_{S}(\nabla _{\overline{\partial}_{j}}\nabla_{\partial_{i}}\partial_{k}))\\
=& D _{\partial_{i}}^{H_{S}(t)}\nabla _{\overline{\partial}_{j}}^{H_{S}(t)}(T_{S}(\partial_{k}))+\langle Rm (\partial_{i} , \overline{\partial}_{j})\partial_{k} , \overline{\partial}_{s}\rangle_{g}g^{l\bar{s}}T_{S}(\partial_{l})\\
& -F_{H_{S}}(\partial_{i} , \overline{\partial}_{j})\circ T_{S}(\partial_{k}) +T_{S}(\partial_{k})\circ F_{H_{S}}(\partial_{i} , \overline{\partial}_{j})\\
=& D _{\partial_{i}}^{H_{S}(t)}(F_{H_{0, S}}(\partial_{k} , \overline{\partial}_{j})-F_{H_{S}}(\partial_{k} , \overline{\partial}_{j}))\\
& +\langle Rm (\partial_{i} , \overline{\partial}_{j})\partial_{k} , \overline{\partial}_{s}\rangle_{g}g^{l\bar{s}}T_{S}(\partial_{l})\\
& -F_{H_{S}}(\partial_{i} , \overline{\partial}_{j})\circ T_{S}(\partial_{k}) +T_{S}(\partial_{k})\circ F_{H_{S}}(\partial_{i} , \overline{\partial}_{j})\\
=& D _{\partial_{i}}^{H_{S}(t)}(F_{H_{0, S}}(\partial_{k} , \overline{\partial}_{j}))- D _{\partial_{k}}^{H_{S}(t)}(F_{H_{S}}(\partial_{i} , \overline{\partial}_{j}))\\
&+\langle Rm (\partial_{i} , \overline{\partial}_{j})\partial_{k} , \overline{\partial}_{s}\rangle_{g}g^{l\bar{s}}T_{S}(\partial_{l})\\
& -F_{H_{S}}(\partial_{i} , \overline{\partial}_{j})\circ T_{S}(\partial_{k}) +T_{S}(\partial_{k})\circ F_{H_{S}}(\partial_{i} , \overline{\partial}_{j}).\\
\end{split}
\end{equation}
Combining the above equalities and recalling $D _{\partial_{i}}^{H_{S}(t)}-D _{\partial_{i}}^{H_{0, S}}=T_{S}(\partial_{i})$, we get
\begin{equation}\label{TS1}
\begin{split}
\Delta |T_{S}|_{H_{S}}^{2} =& 2|\nabla ^{H_{S}(t)}T_{S} |_{H_{S}(t)}^{2} +2Ric_{\omega}(\partial _{k}, \overline{\partial }_{s})g^{k\bar{i}}g^{l\bar{s}}\tr(T_{S}(\partial_{l})H_{S}^{-1}\overline{(T_{S}(\partial_{i}))^{T}}H_{S})\\
& -2Re\langle[\sqrt{-1}\Lambda_{\omega }F_{H_{S}}, T_{S}], T_{S}\rangle_{H_{S}(t)}\\
&-4Re\langle\partial_{H_{S}}(\sqrt{-1}\Lambda_{\omega }F_{H_{S}}), T_{S}\rangle_{H_{S}(t)}\\
& +4Re \{g^{i\bar{j}}g^{k\bar{l}} \langle[T_{S}(\partial_{i}),F_{H_{0, S}}(\partial_{k} , \overline{\partial}_{j}) ], T_{S}(\partial_{l}) \rangle_{H_{S}(t)}\}\\
&+4Re\langle\partial_{H_{0, S}}(\sqrt{-1}\Lambda_{\omega }F_{H_{0, S}}), T_{S}\rangle_{H_{S}(t)}.\\
\end{split}
\end{equation}
From (\ref{S0}), (\ref{TS0}) and (\ref{TS1}), it follows that
\begin{equation}\label{TS}
\begin{split}
&(\Delta -\frac{\partial }{t})|T_{S}|_{H_{S}}^{2}\\=&  2|\nabla ^{H_{S}(t)}T_{S} |_{H_{S}(t)}^{2} +2Ric_{\omega}(\partial _{k}, \overline{\partial }_{s})g^{k\bar{i}}g^{l\bar{s}}\tr(T_{S}(\partial_{l})H_{S}^{-1}\overline{(T_{S}(\partial_{i}))^{T}}H_{S})\\
& -2Re\langle[\sqrt{-1}\Lambda_{\omega }(\gamma \wedge \gamma^{\ast}), T_{S}], T_{S}\rangle_{H_{S}(t)}\\
&-4Re\langle\partial_{H_{S}}(\sqrt{-1}\Lambda_{\omega }(\gamma \wedge \gamma^{\ast})), T_{S}\rangle_{H_{S}(t)}\\
& +4Re \{g^{i\bar{j}}g^{k\bar{l}} \langle[T_{S}(\partial_{i}),F_{H_{0, S}}(\partial_{k} , \overline{\partial}_{j}) ], T_{S}(\partial_{l}) \rangle_{H_{S}(t)}\}\\
&+4Re\langle\partial_{H_{0, S}}(\sqrt{-1}\Lambda_{\omega }F_{H_{0, S}}), T_{S}\rangle_{H_{S}(t)},\\
\end{split}
\end{equation}
and similarly
\begin{equation}\label{TQ}
\begin{split}
&(\Delta -\frac{\partial }{t})|T_{Q}|_{H_{Q}}^{2}\\=&  2|\nabla ^{H_{Q}(t)}T_{Q} |_{H_{Q}(t)}^{2} +2Ric_{\omega}(\partial _{k}, \overline{\partial }_{s})g^{k\bar{i}}g^{l\bar{s}}\tr(T_{Q}(\partial_{l})H_{Q}^{-1}\overline{(T_{Q}(\partial_{i}))^{T}}H_{Q})\\
& -2Re\langle[\sqrt{-1}\Lambda_{\omega }(\gamma^{\ast} \wedge \gamma ), T_{Q}], T_{Q}\rangle_{H_{Q}(t)}\\
&-4Re\langle\partial_{H_{Q}}(\sqrt{-1}\Lambda_{\omega }(\gamma^{\ast} \wedge \gamma )), T_{Q}\rangle_{H_{Q}(t)}\\
&+4Re \{g^{i\bar{j}}g^{k\bar{l}} \langle[T_{Q}(\partial_{i}),F_{H_{0, Q}}(\partial_{k} , \overline{\partial}_{j}) ], T_{Q}(\partial_{l}) \rangle_{H_{Q}(t)}\}\\
&+4Re\langle\partial_{H_{0, Q}}(\sqrt{-1}\Lambda_{\omega }F_{H_{0, Q}}), T_{Q}\rangle_{H_{Q}(t)}.\\
\end{split}
\end{equation}

\medskip

\begin{proposition} Let $H(t)$ be the solution of the heat flow (\ref{Dflow}) with initial metric $H_{0}$ on $(E, \overline{\partial} _{E})$, let $T_{S}(t)$ and $T_{Q}(t)$ be defined by (\ref{DTS1}) and (\ref{DTQ1}), and $\gamma (t)$ be the second fundamental form. Then for any compact subset $U\subset M\setminus \Sigma_{alg}$, there exists a constant $\tilde{C}_{1, U}$ such that
\begin{equation}\label{C101}
\sup_{(x, t)\in U\times [0, +\infty )} (|T_{S}(t)|_{H_{S}(t)}^{2}+|T_{Q}(t)|_{H_{Q}(t)}^{2}+|\gamma (t)|_{H_{S, Q}(t)}^{2})(x, t) \leq \tilde{C}_{1, U}.
\end{equation}
\end{proposition}

\medskip

{\bf Proof. }
Let $\delta $ be small enough so that $U\subset M\setminus B_{\Sigma }(\delta )$. Since the second fundamental form $\overline{\partial }_{S\otimes Q^{\ast}}\gamma_{0} =0$,  for any point $P\in U$, we have a domain $B_{P}(4r)$ and a local section  $G_{0}\in \Gamma (B_{P}(4r); S\otimes Q^{\ast})$ such that $\gamma_{0}=\overline{\partial }_{S\otimes Q^{\ast}} G_{0}$, where $r< \frac{1}{8}\delta $.
Set
\begin{equation}
\rho_{0}=\left (
\begin{matrix}
\Id_{S} &  -G_{0} \\
0   & \Id_{Q}\\
\end{matrix}
\right ).
\end{equation}
It is easy to check that
\begin{equation}
\begin{split}
&\left (
\begin{matrix}
\Id_{S} &  G_{0} \\
0   & \Id_{Q}\\
\end{matrix}
\right )\left (
\begin{matrix}
\overline{\partial}_{S} &  \gamma _{0} \\
0   & \overline{\partial}_{Q}\\
\end{matrix}
\right )\left (
\begin{matrix}
\Id_{S} &  -G_{0} \\
0   & \Id_{Q}\\
\end{matrix}
\right )\\
=&\left (
\begin{matrix}
\overline{\partial}_{S} &  \gamma _{0}- \overline{\partial }_{S}\circ G_{0}+G_{0}\circ \overline{\partial }_{Q}\\
0   & \overline{\partial}_{Q}\\
\end{matrix}
\right ).
\end{split}
\end{equation}
So there exists a local bundle isomorphism $\rho_{0}: S\oplus Q \rightarrow S\oplus Q$ such that
\begin{equation}
\rho_{0}^{\ast }\left (
\begin{matrix}
\overline{\partial}_{S} &  \gamma _{0} \\
0   & \overline{\partial}_{Q}\\
\end{matrix}
\right )=\left (
\begin{matrix}
\overline{\partial}_{S} &  0 \\
0   & \overline{\partial}_{Q}\\
\end{matrix}
\right ),
\end{equation}
and locally
\begin{equation}
(f_{H_{0}}\circ \rho_{0})^{\ast }(\overline{\partial }_{E})=\left (
\begin{matrix}
\overline{\partial}_{S} &  0 \\
0   & \overline{\partial}_{Q}\\
\end{matrix}
\right ).
\end{equation}

Define a local Hermitian metric $\overline{H}_{0, E}$ on $E$ by
\begin{equation}
\overline{H}_{0, E}\doteq (\rho_{0}^{-1}\circ f_{H_{0}}^{-1})^{\ast }(H_{0, S}\oplus H_{0, Q}),
\end{equation}
and set
\begin{equation}
T_{E}(t)=D_{(H(t), \overline{\partial}_{E})}-D_{(\overline{H}_{0, E}, \overline{\partial}_{E})}.
\end{equation}
Now we get
\begin{equation}
(f_{H_{0}}\circ \rho_{0})^{\ast }D_{(\overline{H}_{0, E}, \overline{\partial }_{E})}=\left (
\begin{matrix}
\overline{\partial}_{S} +\partial_{H_{0, S}}&  0 \\
0   & \overline{\partial}_{Q}+\partial_{H_{0, Q}}\\
\end{matrix}
\right ),
\end{equation}
\begin{equation}
(\rho_{0}^{-1}\circ f_{H_{0}}^{-1}\circ f_{H(t)})=\left (
\begin{matrix}
\Id_{S} &  G+G_{0} \\
0   & \Id_{Q}\\
\end{matrix}
\right ),
\end{equation}
\begin{equation}
\begin{split}
&(\rho_{0}^{-1}\circ f_{H_{0}}^{-1}\circ f_{H(t)})^{\ast }\left (
\begin{matrix}
\overline{\partial}_{S} +\partial_{H_{0, S}}&  0 \\
0   & \overline{\partial}_{Q}+\partial_{H_{0, Q}}\\
\end{matrix}
\right )\\
=&\left (
\begin{matrix}
\Id_{S} &  -G-G_{0} \\
0   & \Id_{Q}\\
\end{matrix}
\right )\left (
\begin{matrix}
\overline{\partial}_{S} +\partial_{H_{0, S}}&  0 \\
0   & \overline{\partial}_{Q}+\partial_{H_{0, Q}}\\
\end{matrix}
\right )\left (
\begin{matrix}
\Id_{S} &  G+G_{0} \\
0   & \Id_{Q}\\
\end{matrix}
\right )\\
=&\left (
\begin{matrix}
\overline{\partial}_{S} +\partial_{H_{0, S}}&  (\overline{\partial}_{S\otimes Q^{\ast}}+\partial_{(H_{0,S}, H_{0, Q})})(G+G_{0}) \\
0   & \overline{\partial}_{Q}+\partial_{H_{0, Q}}\\
\end{matrix}
\right )\\
=&\left (
\begin{matrix}
\overline{\partial}_{S} +\partial_{H_{0, S}}&  \gamma (t)+\partial_{(H_{0, S}, H_{0, Q})}(G+G_{0}) \\
0   & \overline{\partial}_{Q}+\partial_{H_{0, Q}}\\
\end{matrix}
\right ),
\end{split}
\end{equation}
and
\begin{equation}
\begin{split}
&f_{H(t)}^{\ast}(T_{E}(t))\\
=&
f_{H(t)}^{\ast}(D_{(H(t), \overline{\partial}_{E})}-D_{(\overline{H}_{0, E}, \overline{\partial}_{E})})\\
=&f_{H(t)}^{\ast}(D_{(H(t), \overline{\partial}_{E})})-(\rho_{0}^{-1}\circ f_{H_{0}}^{-1}\circ f_{H(t)})^{\ast }\circ (f_{H_{0}}\circ \rho_{0})^{\ast }(D_{(\overline{H}_{0, E}, \overline{\partial}_{E})})\\
=&\left (
\begin{matrix}
\overline{\partial}_{S} +\partial_{H_{S}(t)}&  \gamma (t) \\
-\gamma^{\ast} (t)   & \overline{\partial}_{Q}+\partial_{H_{Q}(t)}\\
\end{matrix}
\right )-\left (
\begin{matrix}
\overline{\partial}_{S} +\partial_{H_{0, S}}&  \gamma (t)+\partial_{(H_{0, S}, H_{0, Q})}(G+G_{0}) \\
0   & \overline{\partial}_{Q}+\partial_{H_{0, Q}}\\
\end{matrix}
\right )\\
=&\left (
\begin{matrix}
T_{S}(t) &  -\partial_{(H_{0, S}, H_{0, Q})}(G+G_{0}) \\
-\gamma^{\ast} (t)  & T_{Q}(t)\\
\end{matrix}
\right ),
\end{split}
\end{equation}
where $\partial $ denotes the $(1, 0)$ part of the Chern connection, and $T_{S}=\partial_{H_{S}(t)}-\partial_{H_{0, S}}$, $T_{Q}=\partial_{H_{Q}(t)}-\partial_{H_{0, Q}}$. Then it holds that
\begin{equation}\label{TE0}
\begin{split}
|T_{E}(t)|_{H(t)}^{2}=&|f_{H(t)}^{\ast}(T_{E}(t))|_{f_{H(t)}^{\ast}(H(t))}^{2}\\
=&\Big|\left (
\begin{matrix}
T_{S}(t) &  -\partial_{(H_{0, S}, H_{0, Q})}(G+G_{0}) \\
-\gamma^{\ast} (t)  & T_{Q}(t)\\
\end{matrix}
\right )\Big|_{f_{H(t)}^{\ast}(H(t))}^{2}\\
=&|T_{S}(t)|_{H_{S}(t)}^{2}+|T_{Q}(t)|_{H_{Q}(t)}^{2}\\
&+|\gamma (t)|_{H_{S, Q}(t)}^{2}+|\partial_{(H_{0, S}, H_{0, Q})}(G+G_{0})|_{H_{S, Q}(t)}^{2}.
\end{split}
\end{equation}

On the other hand, calculating in the same way as that in (\ref{TS}) , we have the following local parabolic inequality
\begin{equation}\label{TE}
\begin{split}
&(\Delta -\frac{\partial }{t})|T_{E}|_{H_{E}}^{2}\\=&  2|\nabla ^{H_{E}(t)}T_{E} |_{H_{E}(t)}^{2} +2Ric_{\omega}(\partial _{k}, \overline{\partial }_{s})g^{k\bar{i}}g^{l\bar{s}}\tr(T_{E}(\partial_{l})H_{E}^{-1}\overline{(T_{E}(\partial_{i}))^{T}}H_{E})\\
& +4Re \{g^{i\bar{j}}g^{k\bar{l}} \langle[T_{E}(\partial_{i}),F_{\overline{H}_{0, E}}(\partial_{k} , \overline{\partial}_{j}) ], T_{E}(\partial_{l}) \rangle_{H_{E}(t)}\}\\
&+4Re\langle\partial_{K_{E}}(\sqrt{-1}\Lambda_{\omega }F_{\overline{H}_{0, E}}), T_{E}\rangle_{H_{E}(t)},\\
\end{split}
\end{equation}
and then
\begin{equation}\label{TE1}
\begin{split}
& (\Delta -\frac{\partial }{\partial t}) |T_{E}|_{H_{E}(t)}^{2}\\
\geq & 2 |\nabla^{H_{E}(t)}T_{E}|_{H_{E}(t)}^{2}-c_{1}|T_{E}|_{H_{E}(t)}^{2}-c_{2}\\
\end{split}
\end{equation}
on $B_{P}(2r)$, where constants $c_{1}$ and $c_{2}$ depend only on $\max_{\overline{B_{P}(2r)}}|F_{\overline{H}_{0, E}}|$ and the lower bound of Ricci curvature of $(M, \omega )$.
Let $\varphi _{1}$, $\varphi _{2}$ be nonnegative cut-off functions satisfying
\begin{equation}
 \varphi _{1}(x)=\left \{\begin{split} &1,   \quad  x\in B_{P}(\frac{r}{4}),\\
&0,  \quad x\in M\setminus B_{P}(\frac{r}{2});\\
\end{split}
\right.
\end{equation}
\begin{equation}
 \varphi _{2}(x)=\left \{\begin{split} &1,   \quad  x\in B_{P}(\frac{r}{2}),\\
&0,  \quad x\in M\setminus B_{P}(r);\\
\end{split}
\right.
\end{equation}
and \begin{equation}|d\varphi _{i}|_{\omega}^{2},\ |\Delta \varphi _{i}|\leq\frac{C}{r^{2}} \end{equation} for $i=1, 2$.
We consider the following test function:
\begin{equation}
\zeta_{1} (\cdot , t)=\varphi_{1}^{2}|T_{E}|_{H_{E}(t)}^{2}+a\varphi_{2}^{2}(|G+G_{0}|_{H(t)}^{2}+\tr \hat{h}_{S}(t)+\tr \hat{h}_{Q}^{-1}(t) ),
\end{equation}
where constant $a$ will be chosen large enough.

Let $\zeta_{1} (q, t_{0})=\max_{M\times [0, t_{1}]}\zeta_{1}$, by the definition of $\varphi_{i}$ and the locally uniform estimates of $|G+G_{0}|_{H(t)}^{2}$, $\tr \hat{h}_{S}(t)$, $\tr \hat{h}_{Q}^{-1}(t)$ (i.e. (\ref{CG2}), (\ref{G151}) and (\ref{C0}) ),  we can suppose that
\[
(q, t_{0}) \in B_{P}(\frac{r}{2})\times (0, t_{1}].
\]
By a similar argument in  (\ref{TTTT1}) and (\ref{QQQQ1}), and the definitions of $T_{S}$ and $T_{Q}$, we have
\begin{equation}
(\Delta -\frac{\partial }{\partial t})\tr \hat{h}_{S}\geq c_{3}|T_{S}(t)|_{H_{S}(t)}^{2} -\tilde{c}_{3}
\end{equation}
and
\begin{equation}
(\Delta -\frac{\partial }{\partial t})\tr \hat{h}_{Q}^{-1}\geq c_{4}|T_{Q}(t)|_{H_{Q}(t)}^{2}-\tilde{c}_{4}
\end{equation}
on $B_{P}(2r)$, where $c_{3}$, $\tilde{c}_{3}$, $c_{4}$ and $\tilde{c}_{4}$ are constants depending only on $\max_{\overline{B_{P}(2r)}}(|F_{H_{0, S}}|+|F_{H_{0, Q}}|)$, the local $C^{0}$ bound of $\hat{h}_{S}$ and $\hat{h}_{Q}$.
Let $\tilde{c}_{1}=\sup_{B_{P}(2r)\times [0, +\infty )}|G(t)+G_{0}|_{H(t)}^{2}$. Choosing $a\geq \frac{2c_{1}+8\tilde{c}_{1}+5Cr^{-2}+1}{\min \{c_{3}, c_{4}, 2\}}$, and using (\ref{G6}), (\ref{TE0}), (\ref{TE1}),
 at the maximum point $(q, t_{0})$, we have
\begin{equation}
0\geq  (\Delta -\frac{\partial }{\partial t}) \zeta_{1} \geq  |T_{E}|_{H_{E}(t)}^{2}-c_{5},\\
\end{equation}
where constant $c_{5}$ depends only on the local $C^{0}$ bound of $\tilde{h}_{S}$ and $\tilde{h}_{Q}$, $r^{-2}$, $G_{0}$, $\max_{\overline{B_{P}(2r)}}(|F_{H_{0, S}}|+|F_{H_{0, Q}}|)$ and the lower bound of Ricci curvature of $(M, \omega )$.
We get a uniform constant $c_{6}$ such that
\begin{equation}
\sup_{(x, t)\in M\times [0, +\infty )}\zeta_{1} (x, t)\leq c_{6},
\end{equation}
and then
\begin{equation}
\begin{split}
&\sup_{x\in B_{P}(\frac{r}{4})\times [0, +\infty ) } (|T_{S}(t)|_{H_{S}(t)}^{2}+|T_{Q}(t)|_{H_{Q}(t)}^{2}+|\gamma (t)|_{H_{S, Q}(t)}^{2})(x, t)\\
\leq & \sup_{x\in B_{P}(\frac{r}{4})\times [0, +\infty )} |T_{E}|_{H(t)}(x, t)\leq c_{6}.\\
\end{split}
\end{equation}
 Since $U$ is compact, by finite covering $\{B_{P_{i}}(r_{i})\}_{i=1}^{N}$, there is a constant $c_{7}$ such that
\begin{equation}
\sup_{(x, t)\in U\times [0, +\infty )} (|T_{S}(t)|_{H_{S}(t)}^{2}+|T_{Q}(t)|_{H_{Q}(t)}^{2}+|\gamma (t)|_{H_{S, Q}(t)}^{2})(x, t)\leq c_{7}.
\end{equation}

\hfill $\Box$ \\

\section{Local Curvature estimate }
\setcounter{equation}{0}

From (\ref{r3}), (\ref{TS}) and (\ref{TQ}), we see
\begin{equation}\label{C1}
\begin{split}
&(\Delta -\frac{\partial }{\partial t}) (|\gamma |_{H(t)}^{2}+|T_{S}|_{H_{S}(t)}^{2}+|T_{Q}|_{H_{Q}(t)}^{2})\\
\geq &\frac{1}{4}(|F_{H_{S}}|_{H_{S}(t)}^{2}+|F_{H_{Q}}|_{H_{Q}(t)}^{2}+2|\partial _{H(t)} \gamma |_{H(t)}^{2})\\
& -\hat{C}_{1}|\gamma |_{H(t)}^{2}(|\gamma |_{H(t)}^{2}+|T_{S}|_{H_{S}(t)}^{2}+|T_{Q}|_{H_{Q}(t)}^{2})\\
& -\hat{C}_{2}(|\gamma |_{H(t)}^{2}+|T_{S}|_{H_{S}(t)}^{2}+|T_{Q}|_{H_{S}(t)}^{2})-\hat{C}_{3}\\
\end{split}
\end{equation}
on $M\setminus B_{\Sigma }(\delta )$, where constants $\hat{C}_{i}$  depend only on the uniform local $C^{0}$ bound of $\hat{h}_{S}$ and $\hat{h}_{Q}$, $H_{0, S}$, $H_{0, Q}$ and the lower bound of Ricci curvature of $(M, \omega )$.
Let's recall
\begin{equation}\label{Fi}
|F_{H(t)}|_{H(t)}^{2}=|F_{H_{S}(t)}-\gamma \wedge \gamma^{\ast}|_{H_{S}(t)}^{2}+|F_{H_{Q}(t)}-\gamma^{\ast} \wedge \gamma |_{H_{Q}(t)}^{2}+2|\partial_{H(t)}\gamma |_{H(t)}^{2},
\end{equation}
and
\begin{equation}\label{FFFFF}
\begin{split}
&(\Delta -\frac{\partial }{\partial t})|F_{H(t)}|_{H(t)}^{2}\\
\geq &2|\nabla_{H(t)}F_{H(t)}|_{H(t)}^{2}-\hat{C}_{4}(|F_{H(t)}|_{H(t)} + |Rm(\omega)|_{\omega})|F_{H(t)}|_{H(t)}^{2},
\end{split}
\end{equation}
where constant $\hat{C}_{4}$ depends only on the dimension of $M$. The proof of the inequality (\ref{FFFFF}) can be found in Siu's lectures (page 31 in \cite{Siu01}).

\medskip

{\bf Proof of Theorem \ref{UCE}}
For simplicity, we denote
\begin{equation}
\nu =|\gamma |_{H(t)}^{2}+|T_{S}|_{H_{S}(t)}^{2}+|T_{Q}|_{H_{Q}(t)}^{2}.
\end{equation}
By the estimate (\ref{C101}), we can choose a constant $\hat{C}_{5}$ such that
\begin{equation}
0<\frac{1}{2}\hat{C}_{5}\leq \hat{C}_{5}-\nu (x, t)\leq \hat{C}_{5}
\end{equation}
for all $(x, t)\in M\setminus B_{\Sigma }(\delta )\times [0, +\infty )$. Let $\varphi _{3}$ and $\varphi _{4}$ be nonnegative cut-off functions satisfying
\begin{equation}
 \varphi _{3}(x)=\left \{\begin{split} &1,   \quad  x\in M\setminus B_{\Sigma}(4\delta ),\\
&0,  \quad x\in  B_{\Sigma}(2\delta );\\
\end{split}
\right.
\end{equation}
\begin{equation}
 \varphi _{4}(x)=\left \{\begin{split} &1,   \quad  x\in M\setminus B_{\Sigma}(2\delta ),\\
&0,  \quad x\in  B_{\Sigma }(\delta);\\
\end{split}
\right.
\end{equation}
and \begin{equation}|d\varphi _{i}|_{\omega}^{2},\ |\Delta \varphi _{i}|\leq\frac{C}{r^{2}} \end{equation} for all $i=3, 4$.

Now we consider the following test function
\begin{equation}
\zeta_{2} =\varphi_{3}^{2}\frac{|F_{H(t)}|_{H(t)}^{2}}{\hat{C}_{5}-\nu }+b\varphi_{4}\nu .
\end{equation}

Let $\zeta_{2} (P, t_{0})=\max_{M\times [0, t_{1}]}\zeta$, by the definition of the cut-off functions $\varphi_{i}$ and the local uniform estimate (\ref{C101}), we can suppose that
\[
(P, t_{0})\in M\setminus B_{\Sigma }(2\delta )\times (0, t_{1}].
\]
At the maximum point $(P, t_{0})$, we have
\begin{equation}\label{Fii}
\begin{split}
&(\Delta -\frac{\partial }{\partial t})\zeta_{2} \\
= & \varphi_{3}^{2}\frac{1}{\hat{C}_{5}-\nu }(\Delta -\frac{\partial }{\partial t})|F_{H(t)}|_{H(t)}^{2}+\varphi_{3}^{2}\frac{|F_{H(t)}|_{H(t)}^{2}}{(\hat{C}_{5}-\nu )^{2} }(\Delta -\frac{\partial }{\partial t})\nu \\
 & -\varphi_{3}^{2}\frac{2}{\hat{C}_{5}-\nu }\nabla (\frac{|F_{H(t)}|_{H(t)}^{2}}{\hat{C}_{5}-\nu })\cdot \nabla (\hat{C}_{5}-\nu )+b (\Delta -\frac{\partial }{\partial t})\nu \\
 & +\Delta \varphi_{3}^{2} \frac{|F_{H(t)}|_{H(t)}^{2}}{\hat{C}_{5}-\nu } +2 \nabla \varphi_{3}^{2} \cdot \nabla \frac{|F_{H(t)}|_{H(t)}^{2}}{\hat{C}_{5}-\nu }\\
\end{split}
\end{equation}
and
\begin{equation}\label{Fii2}
\nabla (\varphi_{3}^{2}\frac{|F_{H(t)}|_{H(t)}^{2}}{\hat{C}_{5}-\nu })+b\nabla \nu=0.
\end{equation}
Putting (\ref{Fii2}) into (\ref{Fii}), choosing the constants $\hat{C}_{5}$ and $b$ large enough and using the formulas (\ref{C1}), (\ref{Fi}), and (\ref{FFFFF}), at the  maximum point $(p, t_{0})$, we have
\begin{equation}\label{Fii3}
\begin{split}
&(\Delta -\frac{\partial }{\partial t})\zeta_{2} \\
= & \varphi_{3}^{2}\frac{1}{\hat{C}_{5}-\nu }(\Delta -\frac{\partial }{\partial t})|F_{H(t)}|_{H(t)}^{2}+\varphi_{3}^{2}\frac{|F_{H(t)}|_{H(t)}^{2}}{(\hat{C}_{5}-\nu )^{2} }(\Delta -\frac{\partial }{\partial t})\nu \\
 & -\frac{2b}{\hat{C}_{5}-\nu } |\nabla \nu |^{2} +b (\Delta -\frac{\partial }{\partial t})\nu \\
 & +\Delta \varphi_{3}^{2} \frac{|F_{H(t)}|_{H(t)}^{2}}{\hat{C}_{5}-\nu } + \frac{2}{\hat{C}_{5}-\nu } \nabla \varphi_{3}^{2} \cdot \nabla |F_{H(t)}|_{H(t)}^{2}\\
\geq & |F_{H(t)}|_{H(t)}^{2} -\hat{C}_{6},\\
\end{split}
\end{equation}
where $\hat{C}_{6}$ is a positive constant depending only on the local uniform bound of $\nu $, $\delta^{-1}$ and the curvature of $(M, \omega )$.
So we obtain
\begin{equation}
|F_{H(t)}|_{H(t)}^{2}(P, t_{0})\leq \hat{C}_{6}.
\end{equation}
Then there is a constant $\hat{C}_{7}$ such that
\begin{equation}
\sup_{M\setminus B_{\Sigma }(4\delta )\times [0, +\infty )}|F_{H(t)}|_{H(t)}^{2}\leq \hat{C}_{7}.
\end{equation}
This completes the proof of Theorem \ref{UCE}.

\hfill $\Box$ \\

\hspace{0.4cm}

\hspace{0.3cm}

\end{document}